\documentclass[a4paper]{article}
\pdfoutput=1
\usepackage{amsmath,amsfonts,amsthm}
\usepackage{graphicx} 
\usepackage[unicode,final]{hyperref} 
\usepackage{cite}
\usepackage{caption}
\usepackage{mathtools}
\usepackage{newtxtext}
\usepackage[varvw]{newtxmath}
\usepackage{bm}
\usepackage[T1]{fontenc}
\usepackage[utf8]{inputenc}
\usepackage{color}

\newcommand{\fz}{\frac}

\newcommand{\lA}{\langle}
\newcommand{\rA}{\rangle}
\newcommand{\rarrow}{\rightarrow}

\newtheorem{Hyp}{Hypothesis}
\newtheorem{Def}{Definition}
\newtheorem{Lem}{Lemma}
\newtheorem{Prop}{Proposition}
\newtheorem{Thm}{Theorem}
\newtheorem{Cor}{Corollary}
\newtheorem{Rmk}{Remark}
\newtheorem{Ex}{Example}
\allowdisplaybreaks[4]
\newcommand{\R}{\mathbb{R}}
\newcommand{\N}{\mathbb{N}}

\newcommand{\defeq}{\coloneqq}
\newcommand{\eqdef}{\eqqcolon}
\providecommand{\keywords}[1]{\textit{Keywords:} #1}
\providecommand{\msc}[1]{\textit{2010 MSC:} #1}

\makeatletter

\renewcommand{\proofname}{Proof}
\makeatother

\title{Error Estimates for First- and Second-Order Lagrange--Galerkin Moving Mesh Schemes for the One-Dimensional Convection-Diffusion Equation}

\author{Kharisma Surya Putri\footnote{Division of Mathematical and Physical Sciences, Kanazawa University, {\tt kharismasp@stu.kanazawa-u.ac.jp}} \and
Tatsuki Mizuochi\footnote{Division of Mathematical and Physical Sciences, Kanazawa University, {\tt green9014@stu.kanazawa-u.ac.jp}} \and 
Niklas Kolbe\footnote{Institute of Geometry and Practical Mathematics, RWTH Aachen University, {\tt kolbe@igpm.rwth-aachen.de}} \and
Hirofumi Notsu\footnote{Faculty of Mathematics and Physics, Kanazawa University, {\tt notsu@se.kanazawa-u.ac.jp}}
}

\date{}
\begin{document}
\maketitle
\begin{abstract}
 A new moving mesh scheme based on the Lagrange--Galerkin method for the approximation of the one-dimensional convection-diffusion equation is studied. The mesh movement, which is prescribed by a discretized dynamical system for the nodal points, follows the direction of convection. It is shown that under a restriction of the time increment the mesh movement cannot lead to an overlap of the elements and therefore an invalid mesh. For the linear element, optimal error estimates in the $\ell^\infty(L^2) \cap \ell^2(H_0^1)$ norm are proved in case of both, a first-order backward Euler method and a second-order two-step method in time. These results are based on new estimates of the time dependent interpolation operator derived in this work. Preservation of the total mass is verified for both choices of the time discretization. Numerical experiments are presented that confirm the error estimates and demonstrate that the proposed moving mesh scheme can circumvent limitations that the Lagrange--Galerkin method on a fixed mesh exhibits. \\
\medskip\\
\\
\msc{65M15, 65M25, 65M50}
\end{abstract}

\section{INTRODUCTION}
Apart from their application in fluid flow problems convection-diffusion equations have recently been extensively used in the modeling of chemical and biological processes such as pollutant transport, immune system dynamics and cancer growth, see e.g. \cite{anderson2000mathemmodeltumourinvasmetas, mcrae1982numer}.
In many of these applications, migration of cells or transport of solutes in fluids are involved, which leads to convection-dominant problems. In these cases the Galerkin finite element scheme can produce oscillating solutions. Hence, a plethora of extended and alternative numerical methods have been developed to perform stable computation, e.g., upwind methods \cite{baba1981,hughes1987errat,tabata1977} and characteristics methods \cite{colera2021lagrangaler,douglas1982numer,pironneau1989finitelemenmethodfluid,tabata2016lagrangaler}. Among the latter ones the Lagrange--Galerkin (LG) method has been shown to be an effective and efficient method to deal with convection-dominated problems, see for example \cite{benitez2012numerlagrangaler,benitez2012IInumerlagrangaler,bermejo2012modiflagrangaler,ewing1981multisgalermethod,futai2022masspresertwo,rui2002,rui2010massconsercharac}. Not only the convection-dominant nature but also the rich dynamics of the biological applications, which include traveling waves and aggregation phenomena, pose a challenge for the numerical schemes. Different approaches related to mesh adaptation have been recently proposed to improve the accuracy in these cases, including a mass-transport approach for the one-dimensional problem~\cite{carrillo2019kellersegel} and adaptive mesh refinement~\cite{acharya1990, kolbe2022adaptrectanmesh}.

The use of the LG method in convection-diffusion problems has the advantage that time increments do not need to be restricted by any CFL-type condition, which has also motivated applications to the Navier--Stokes equations and models of viscoelastic fluid flow to name only a few \cite{notsu2016errorestimstabil,lukacova-medvidova2017numeroseenpeter, lukacova-medvidova2017numeroseenpeterII}. Various LG methods have been proposed for convection-diffusion problems: Some have been concerned with higher order approximations in time using single- and multi-step methods, e.g., \cite{benitez2012numerlagrangaler,benitez2012IInumerlagrangaler,bermejo2015modiflagrangaler}. Others have focused on deriving schemes that maintain the mass balance on the discrete level, e.g. \cite{colera2020lagrangaler,rui2010massconsercharac}. In~\cite{futai2022masspresertwo,rui2010massconsercharac} such mass-preserving LG schemes of first- and second-order in time have been proposed and error estimates have been provided. As the LG scheme relies on an upwind-interpolation of the numerical solution that follows the velocity field backwards in time a promising approach is to introduce mesh movement along the velocity field. From a computational point of view this might ease the identification of the upwind points and reduce the interpolation error.

Different kinds of moving mesh methods have been considered to numerically solve convection-diffusion problems. A common approach is the redistribution of mesh cells according to a monitor function that depends on local features of the numerical solution or an a posteriori error estimator, see e.g.~\cite{acharya1990}. In other approaches separate moving mesh PDEs and transformations obtained from the solution of the Monge--Kantorovich problem are used, see \cite{huang2011adapt,sulman2011optim}. In the context of hyperbolic balance laws and fluid dynamics a variety of schemes, which also entail mesh movement, has been derived from the Lagrangian formulation of the problem, such as the hydrodynamic \emph{GLACE} scheme \cite{MR2537850} and the arbitrary Lagrangian--Eulerian finite element method \cite{TAKASHI1992115, MR3030556}. In this context high order of accuracy has been achieved by adopting high order essentially non-oscillatory reconstructions, see e.g., \cite{MR2442406, MR3079104}. In some applications error estimates have been derived taking the movement of the mesh into account, e.g. \cite{chrysafinos2008lagran,gelinas1981}. We consider here a Lagrangian grid approach based on a dynamical system for the nodal points that both exhibits significant benefits over static grids in numerical simulations and allows for an error analysis of the extended LG scheme. 

In this work we are concerned with a moving mesh approach within LG schemes of first- and second-order in time and corresponding error estimates. We introduce the Lagrange--Galerkin Moving Mesh (LGMM) schemes, which combine the LG schemes derived in~\cite{futai2022masspresertwo,rui2010massconsercharac} with a moving mesh, in order to improve the performance and efficiency over the original LG schemes with fixed mesh. The mesh movement we consider is inspired by \cite{carrillo2019kellersegel}; although we focus on the one-dimensional case in this work our mesh movement is expressed in a form suitable for higher dimensional cases. We derive a condition under which the mesh movement is applicable. We generalize the mass-conservation property and stability results from the static grid versions of the schemes to the new LGMM schemes. Based on the idea of temporal derivatives on deforming grids in \cite{jimack1991tempor} we derive bounds for the time derivative of the dynamic interpolation operator, which then allow us to prove optimal error estimates for the LGMM scheme on linear elements in the $\ell^\infty(L^2) \cap \ell^2(H_0^1)$ norm. Moreover, we present numerical experiments that verify the error estimates. They further show that in case of aggregation the LGMM method eliminates oscillations of the numerical solution that the LG method produces.

The rest of the paper is organized as follows: In Section 2 we present the mass-preserving LGMM schemes of first- and second-order in time for convection-diffusion problem. In Section 3, we provide the main results concerned with the mesh movement, the mass-conservation property, the stability, and the error estimates, which are afterwards proven in Section 4. To show the advantages of the LGMM schemes, two numerical simulations are given in Section 5, followed by the conclusions in Section 6. 

%\\\\\\\\\\\\\\\\\\\\\\\\\\\\\\\\\\\\\\\\\\\\\\\\\\\\\\\\\\\\\\\\\\\
\section{LAGRANGE--GALERKIN AND MOVING MESH METHOD}
%\\\\\\\\\\\\\\\\\\\\\\\\\\\\\\\\\\\\\\\\\\\\\\\\\\\\\\\\\\\\\\\\\\\
\subsection{Statement of the Problem}
Let $\Omega=(a,b)$ be a bounded interval in $\R$. We denote by $\Gamma \defeq \partial\Omega$ the two point boundary of $\Omega$ and by $T$ a positive constant. In this paper we use the Lebesgue spaces $L^2(\Omega),L^{\infty}(\Omega)$ and the Sobolev spaces $W^{m,p}(\Omega),W^{1,\infty}_0(\Omega),H^m(\Omega),H^1_0(\Omega)$, for $m\in\mathbb{N}\cup\{0\}$ and $p\in[1,\infty]$. We use the notation $(\cdot,\cdot)$ to represent the $L^2(\Omega)$ inner products for both scalar and vector-valued functions. The norm in $L^2(\Omega)$ is simply denoted as $ \| \cdot \| \defeq \| \cdot \| _{L^2(\Omega)}$. For any normed space $Y$ with norm $ \| \cdot \| _Y$, we define the function spaces $H^m(0,T;Y)$ and $C(0,T;Y)$ consisting of $Y$-valued functions in $H^m(0,T)$ and $C([0,T])$, respectively. For the two real numbers $t_0<t_1$ we introduce the function space 
\begin{equation*}
    Z^m(t_0,t_1)\defeq\{\psi\in C^j(t_0,t_1;H^{m-j}(\Omega));j=0,\dots,m, \| \psi \| _{Z^m(t_0,t_1)}<\infty\},
\end{equation*}
with the norm
\begin{equation*}
 \|  \cdot  \| _{Z^m(t_0,t_1)}\defeq \left(\sum_{j=1}^m  \|  \cdot  \| ^2_{C^j(t_0,t_1;H^{m-j}(\Omega))}\right)^{1/2},
\end{equation*}
and set $Z^m\defeq Z^m(0,T)$. We often omit $\Omega$ and $ [0,T]$ if there is no confusion and write, e.g., $C(L^\infty)$ in place of $C([0,T];L^\infty(\Omega))$. Although we are concerned with a one-dimensional domain we use the general notations $\nabla \defeq \partial_x$, $\nabla \cdot \defeq \partial_x$, $\Delta \defeq \partial_x^2$, and $\frac{\partial}{\partial n} \defeq n \partial_x$ to refer to spatial derivatives in order to allow for a straightforward application of the multi-dimensional theory for LG schemes. We use $c$ and $C$ (with or without subscript or superscript) to denote generic positive constant independent of discretization parameters and solutions.

We consider a convection-diffusion problem, in which we aim to find $\phi:\Omega\times(0,T) \rightarrow\mathbb{R}$ such that
\begin{subequations}\label{prob:c-d}
\begin{alignat}{2}
    \frac{\partial\phi}{\partial t}+\nabla\cdot(u\phi)-\nu\Delta\phi&=f &&\text{in }\Omega\times(0,T),\\
    \nu\frac{\partial\phi}{\partial n}-\phi u \cdot n&=g &&\text{on }\Gamma\times(0,T),\\
    \phi&=\phi^0 &&\text{ in }\Omega,\text{at }t=0,
\end{alignat}
\end{subequations}
where $u:\Omega\times(0,T)\rightarrow\mathbb{R}$, $f:\Omega\times(0,T)\rightarrow\mathbb{R}$, $g:\Gamma\times(0,T)\rightarrow\mathbb{R}$ and $\phi^0:\Omega\rightarrow\mathbb{R}$ are given functions, $n:\Gamma \to \{ \pm 1\}$ is the outward unit normal vector and $\nu >0$ is the diffusion coefficient. 

Let $\Psi\defeq H^1(\Omega)$ and $\Psi'$ be the dual space of $\Psi$. The weak formulation corresponding to problem~\eqref{prob:c-d} is to find $\{\phi(t)=\phi(\cdot,t)\in\Psi; t\in(0,T)\}$ such that for $t\in(0,T)$ the variational equality
\begin{equation}\label{weakform}
\left(\frac{\partial\phi}{\partial t}(t),\psi\right)+a_0(\phi(t),\psi)+a_1(\phi(t),\psi;u(t))=\langle F(t),\psi \rangle, \quad
\forall\psi\in\Psi
\end{equation}
holds in addition to $\phi(0)=\phi^0$. The bilinear forms $a_0(\cdot,\cdot)$ and $a_1(\cdot,\cdot)=a_1(\cdot,\cdot\,;u)$ and for $t \in(0,T)$ in~\eqref{weakform} are defined by
\begin{align*}
    a_0(\phi,\psi)  \defeq \nu(\nabla\phi,\nabla\psi),\qquad
                     a_1(\phi,\psi;u)  \defeq -(\phi,u \nabla\psi)
\end{align*}
and $F(t)\in \Psi'$ is a functional given by
\begin{align*}
\langle F(t),\psi \rangle=(f(t),\psi)+[g(t),\psi]_{\Gamma}, \qquad
    [g(t),\psi]_{\Gamma}  \defeq \int_{\Gamma}g(t)\psi ds,
\end{align*}
for $f(t)=f(\cdot,t)\in L^2(\Omega)$ and $g(t)=g(\cdot,t)\in L^2(\Gamma)$. 

By substituting $1\in\Psi$ into $\psi$ in \eqref{weakform}, and integrating over $(0,t)$, we derive the mass balance identity
\begin{align}\label{mass_id}
\int_{\Omega}\phi(x,t) \, dx = \int_{\Omega}\phi^0(x) \, dx+\int_{0}^{t} \int_{\Omega}f(x,\tau) dx \, d\tau
+\int_{0}^{t}\int_{\Gamma}g(s,\tau) ds \, d\tau
\end{align}
that holds for all $t\in(0,T)$.
This property is desired to be maintained also at the discrete level, which is indeed achieved by the Lagrange--Galerkin schemes of first- and second-order in time proposed in \cite{rui2010massconsercharac} and \cite{futai2022masspresertwo}, respectively.
%\\\\\\\\\\\\\\\\\\\\\\\\\\\\\\\\\\\\\\\\\\\\\\\\\\\\\\\\\\\\\\\\\\\
\subsection{The First-Order Lagrange--Galerkin Scheme}\label{sec:1storderLG}
Let $\Delta t > 0$ be a time increment, $t^n \defeq n\Delta t$ for $n\in\mathbb{Z}$ an equidistant discretization of the time domain and $N_T \defeq \lfloor{T/\Delta t}\rfloor$. For a function~$\rho$ defined in $\Omega\times(0,T)$ and $0\leq t^n\leq T$ the function $\rho(\cdot,t^n)$ in $\Omega$ is denoted by $\rho^n$. 
Let $\mathcal{T}_h^n\defeq \{K^n\},$ $n\in\{0,\ldots,N_T\}$ be a time-dependent partition of $\bar{\Omega}$ with $K^n$ representing an element in $\mathcal{T}_h^n$. Let $h\defeq\max_{n=0,\dots,N_T}\max_{K^n\in \mathcal{T}_h^n} \text{diam} (K^n)$ denote the global mesh size. 
Let $\Psi_h^n \subset \Psi$ be a time varying finite element space defined by
\[
  \Psi_h^n \defeq \left\{\psi_h\in C(\bar{\Omega});\ \psi_{h|K^n} \in P_1(K^n),\forall K^n\in\mathcal{T}_h^n\right\},
\]
where $P_1(K^n)$ is the space of linear polynomial functions on $K^n\in\mathcal{T}_h^n$. Details about the spaces $\Psi_h^n$ and $\mathcal{T}_h^n$ are discussed in Section~\ref{sec:movingmesh}, the rest of this section as well as Section~\ref{sec:2ndorderLG} address the first- and second-order schemes assuming $\Psi_h^n$ is known.

For a given velocity $v:\Omega\rightarrow\mathbb{R}$, we define the upwind point of~$x$ with respect to $v$ and $\Delta t$ using the mapping $X_1(v,\Delta t):\Omega\rightarrow\mathbb{R}$,
\begin{equation*}
    X_1(v,\Delta t)(x) \defeq x-v(x)\Delta t.
\end{equation*}
With respect to the velocity $u$ we define the mapping $X_1^n : \Omega\rightarrow\R$ and its corresponding Jacobian $\gamma^n : \Omega\rightarrow\R$ by
\begin{align*}
X_1^n(x) \defeq X_1(u^n,\Delta t)(x)=x-u^n(x)\Delta t, \qquad \gamma^n(x) \defeq \det\left(\frac{\partial X_1^n}{\partial x}(x)\right).
\end{align*}
Suppose an approximate function $\phi_h^0\in\Psi^0_h$ of $\phi^0$ is given. In the first-order Lagrangian moving mesh scheme we look for $\{\phi_h^n\in\Psi_h^n;\ n=1,\dots,N_T\}$ such that for $n = 1,\ldots, N_T$ it holds 
\begin{align} \label{LG}
\left(\frac{\phi^n_h-\phi^{n-1}_h\circ X^n_1 \gamma^n}{\Delta t}, \psi_h \right) + a_0( \phi_h^n, \psi_h ) = \lA F^n, \psi_h \rA, \quad
\forall \psi_h \in \Psi_h^n.
\end{align}
The functional $F^n\in (\Psi_h^n)^\prime$ on the right hand side of \eqref{LG} is defined by
\[
  \lA F^n, \psi_h \rA  \defeq( f^n, \psi_h ) + [g^n, \psi_h]_\Gamma.
  \]
%\\\\\\\\\\\\\\\\\\\\\\\\\\\\\\\\\\\\\\\\\\\\\\\\\\\\\\\\\\\\\\\\\\\
\subsection{The Second-Order Lagrange--Galerkin Scheme}\label{sec:2ndorderLG}
To obtain a higher order discretization in time we define the additional mapping $\tilde{X}_1^n:\Omega\rightarrow\R$ and its Jacobian $\tilde{\gamma}^n:\Omega\rightarrow\R$ by
\begin{align*}
\tilde{X}_1^n(x) \defeq X_1(u^n,2\Delta t)(x)=x-2u^n(x)\Delta t, \qquad 
\tilde{\gamma}^n(x) \defeq \det\left(\frac{\partial \tilde{X}_1^n}{\partial x}(x)\right).
\end{align*}
Suppose an approximation $\phi_h^0\in\Psi^0_h$ of $\phi^0$ is given. Then the second-order Lagrangian moving mesh scheme aims to find $\{\phi_h^n\in\Psi_h^n;\ n=1,\dots,N_T\}$ satisfying 
\begin{subequations}\label{LG1}
\begin{align} 
\left(\frac{\phi^n_h-\phi^{n-1}_h\circ X^n_1 \gamma^n}{\Delta t}, \psi_h \right) + a_0( \phi_h^n, \psi_h ) &= \lA F^n, \psi_h \rA, \quad
\forall \psi_h \in \Psi_h^n,\quad n=1,\\
\left(\frac{3\phi^n_h-4\phi^{n-1}_h\circ X^n_1 \gamma^n+ \phi^{n-2}_h\circ \tilde{X^n_1}\tilde{\gamma^n}}{2 \Delta t}, \psi_h \right) + a_0( \phi_h^n, \psi_h ) &= \lA F^n, \psi_h \rA, \quad
\forall \psi_h \in \Psi_h^n, \quad n\geq2. \label{eq:lg22}
\end{align}
\end{subequations}
In the following, we rewrite scheme~\eqref{LG1} as
\begin{align} \label{LG2}
( \mathcal{A}_{\Delta t} \phi_h^n, \psi_h ) + a_0( \phi_h^n, \psi_h ) = \lA F^n, \psi_h \rA, \ \ 
\forall \psi_h \in \Psi_h^n,
\end{align}
where, for a series~$\{\rho^n\}_{n=0}^{N_T} \subset \Psi$, the function $\mathcal{A}_{\Delta t} \rho^n: \Omega \to \R$ is given by
\begin{align*}
\mathcal{A}_{\Delta t} \rho^n & \defeq
\left\{
\begin{aligned}
& \mathcal{A}_{\Delta t}^{(1)} \rho^n\defeq\frac{1}{\Delta t}\left(\rho^n-\rho^{n-1} \circ X^n_1 \gamma^n\right),
& n & = 1, \\
&  \mathcal{A}_{\Delta t}^{(2)} \rho^n\defeq\fz{1}{2\Delta t} \left( 3\rho^n-4\rho^{n-1} \circ X_1^n \gamma^n + \rho^{n-2} \circ \tilde{X}_1^n \tilde{\gamma}^n \right),
& n & \ge 2,
\end{aligned}
\right.
\end{align*}

We introduce some discrete norms in the following. Let $Y$ be a normed space, $m\in \{0,\dots,N_T\}$ be an integer, and $\{\rho^n\}^{N_t}_{n=0}\subset Y$. We define the norms $ \| \cdot \| _{\ell^\infty_m (Y)}$ and $ \| \cdot \| _{\ell^2_m(Y)}$ by
\begin{equation*}
     \| \rho \| _{\ell^\infty_m (Y)}\defeq\max_{n=m,\dots,N_T}  \| \rho^n \| _{Y}, \quad
     \| \rho \| _{\ell^2_m(Y)}\defeq\left(\Delta t \sum_{n=m}^{N_T}  \| \rho^n \| ^2_{Y}\right)^{1/2}.
  \end{equation*}
  Additionally, we define the following norm over the time varying finite element spaces
\[
  \| \rho \|_{\ell_m^2(\Psi_h^\prime)} \defeq \left(\Delta t \sum_{n=m}^{N_T}  \| \rho^n \| ^2_{(\Psi_h^n)^\prime}\right)^{1/2}.
  \]
If there is no confusion we omit the subscript i.e, $ \| \rho \| _{\ell^\infty_1(Y)} \eqdef  \| \rho \| _{\ell^\infty(Y)}$ and $ \| \rho \| _{\ell^2_1(Y)} \eqdef \| \rho \| _{\ell^2(Y)}$. 

% \\\\\\\\\\\\\\\\\\\\\\\\\\\\\\\\\\\\\\\\\\\\\\\\\\\\\\\\\\\\\\\\\\\
\subsection{Moving Mesh}\label{sec:movingmesh}
In this section the construction and evolution of the partitions $\mathcal{T}_h^n$ is considered. To this end a moving mesh is employed that in this work is defined as follows.
\begin{Def}\label{def:mm} For a given partition $\{t^n:~n=0,\dots,N_T\}$ of the time domain $[0,T]$ a \emph{moving mesh} of $\Omega \times [0, T]$ is a set of points $\{P_i^n:~i=1,\dots,N_p,~n=0,\dots,N_T\} \subset \bar \Omega$ that satisfy the monotonicity condition
\begin{equation}\label{eq:nooverlap}
  a=P^n_1< P^n_2< \cdots < P^n_{N_p}=b, \qquad  n \in \{0,\dots,N_T\}.
\end{equation}
We refer to $N_p \in \N$ as the number of moving mesh points.
\end{Def}

A moving mesh allows us to define the partitions introduced in Section~\ref{sec:1storderLG} more precisely as
\[
  \mathcal{T}_h^n = \{K_{i}^n:~i=1,\dots, N_p - 1 \}, \qquad K_{i}^n \defeq [P^n_i, P^n_{i+1}], \qquad n \in \{0,\dots,N_T\}
  \]
  and therefore determines the nodal points of the finite element spaces $\Psi_h^n$ for $n\in\{0, \dots, N_T\}$. We note that for any fixed $i\in\{1, \dots, N_p\}$ the series $\{P_i^n\}_{n=0}^{N_T}$ can be considered a time discrete trajectory of the moving point $P_i$. For our analysis we define the velocities of the moving mesh points as
  \begin{equation}
     w^n_i\defeq\frac{P^n_i-P^{n-1}_i}{\Delta t}, \qquad i\in\{1,\dots,N_p\},\quad n\in\{1,\dots,N_t\},
  \end{equation}
  which allow us to introduce the time continuous point trajectories
  \begin{equation}\label{eq:movingpoints}
    P_i(t) \defeq P^{n-1}_i+w^n_i(t-t^{n-1}), \qquad i\in\{1,\dots,N_p\}, \quad t\in[t^{n-1},t^n]. 
  \end{equation}
  In addition, for $t\in[0,T]$  we define $w(x,t)$ as an extension of $w^n_i$, by
\begin{equation}
    w(x,t)\defeq\frac{P_{i+1}(t)-x}{P_{i+1}(t)-P_i(t)}w^n_i+\frac{x-P_i(t)}{P_{i+1}(t)-P_i(t)}w^n_{i+1}, \qquad x \in[P_i(t), P_{i+1}(t)], \quad t\in[t^{n-1},t^n].
\end{equation}
Note that $w \in C^0([0,T];W^{1,\infty}_0(\Omega))$. Also the basis functions of the finite element spaces $\Psi_h^n$ can be naturally extended using the trajectories~\eqref{eq:movingpoints}: for $i\in\{1,\dots, N_p\}$ let $\psi_i(\cdot, t)$ denote the unique function on $\Omega$ that is affine linear restricted to the intervals $[P_j(t), P_{j+1}(t)]$ for $j\in\{1, \dots, N_p -1 \}$ and satisfies $\psi_i(P_j(t), t) = \delta_{ij}$. Then clearly $\{ \psi_i(\cdot, t^n):~i=1,\dots,N_p \}$ is a basis of $\Psi_h^n$.

\subsection{Moving Mesh Method}\label{sec:movingmeshmethod}
In this section we propose a \emph{moving mesh method} that is used to obtain a moving mesh in the sense of Definition~\ref{def:mm} and therefore determines the finite element spaces $\Psi_h^n$ as described in Section~\ref{sec:movingmesh}. Suppose that the points $P_1^0$, \dots, $P_{N_P}^0$ are given and the monotonicity condition~\eqref{eq:nooverlap} is satisfied for $n=0$. The method we propose determines the position of the points $P_i^n$ iteratively by employing a time discretization of the dynamical system
\begin{align} \label{MM_c}
\fz{d \tilde P_i}{dt}(t) = u\bigl( \tilde P_i(t), t \bigr) & + \nu_M \sum_{j=1}^{N_p -1}\nabla\psi_i(t)\rvert_{[\tilde P_j(t), \tilde P_{j+1}(t)]}
\end{align}
with initial data $\tilde P_i(0)=P_i^0$ for $i\in\{1, \dots, N_p\}$. Here, the parameter $\nu_M\geq 0$ accounts for regularization of the moving mesh. 

Applying a linearly implicit time discretization to the continuous problem~\eqref{MM_c} gives rise to our moving mesh method: find $\{P_i^n:~i=1,\dots,N_p,~n=0,\dots,N_T\}$ such that for $n=1,\dots,N_T$ it holds
\begin{subequations} \label{MM_d}
\begin{align}
    \frac{P^{n}_i-P^{n-1}_i}{\Delta t} &=u^{n-1}(P^{n-1}_i)+\nu_M \frac{P^{n}_{i+1}-2P^{n}_i+P^{n}_{i-1}}{(P^{n-1}_i-P^{n-1}_{i-1})(P^{n-1}_{i+1}-P^{n-1}_i)}, \qquad i= 2,\dots,N_p-1, \label{mm_d_1} \\
    P^{n}_1&=a,\qquad P^{n}_{N_p}=b, \label{mm_d_2} \\
   \{P_i^0:~i=1,\dots,N_p\}&\subset\bar{\Omega} \text{ given };\qquad a=P^0_1< P^0_2<...< P^0_{N_p}=b. \label{mm_d_3}
\end{align}
\end{subequations}
The method is inspired by \cite{carrillo2019kellersegel} and can be extended to higher dimensions in a straightforward way. 
In the case $\nu_M= 0$, the transition from $P_i^{n-1}$ to $P_i^n$ due to \eqref{MM_d} and the transition from $P_i^{n-1}$ to $X_1(u^{n-1}, \Delta t)(P_i^{n-1})$ describe movements in opposite directions. In particular, if the velocity field $u$ is smooth we have $ P_i^{n-1} \approx X_1^n(P_i^n)$. Hence, a reduction of the computational costs to identify $X_1^n(P_i^n)$ as well as a decrease of the corresponding interpolation error in the scheme are expected. The main idea of the LGMM method is, to combine the LG schemes~\eqref{LG} and~\eqref{LG2} with the moving mesh method \eqref{MM_d}.

\begin{Rmk}
While the moving mesh method~\eqref{MM_d} leads to a well defined set of nodal points $\{P_i^n:~i=1,\dots,N_p,~n=0,\dots,N_T\}$ it is not clear whether they constitute a moving mesh in the sense of Definition~\ref{def:mm} since the condition~\eqref{eq:nooverlap} might not be satisfied.
\end{Rmk}

\begin{Rmk}
Generally, the coefficient matrix resulting from~\eqref{mm_d_1} is not symmetric.
\end{Rmk}

In fact we show that the method \eqref{MM_d} results in a moving mesh for a suitable choice of $\Delta t$, see Section~\ref{sec:resultsmm}. The other theoretical results we show in this work assume a given moving mesh.
While it is important that a positive distance between neighbor points is maintained in the moving mesh, it needs to be also verified that this distance does not become too large. In particular, with respect to the error estimates that we present in the following section we are interested in the situation that the global mesh size $h$ tends to $0$. This can be realized by employing an equidistant mesh of size $h_0$ at the initial time that is iteratively decreased by increasing the number of moving points and ensuring that the distance between neighboring points does not exceed $C h_0$ over all time instances for a fixed constant $C>0$. In practice, positive $\nu_M$ in scheme \eqref{MM_d} has resulted in a control over the maximal point distance.
% \\\\\\\\\\\\\\\\\\\\\\\\\\\\\\\\\\\\\\\\\\\\\\\\\\\\\\\\\\\\\\\\\\\
\section{MAIN RESULTS}\label{sec:results}
%\\\\\\\\\\\\\\\\\\\\\\\\\\\\\\\\\\\\\\\\\\\\\\\\\\\\\\\\\\\\\\\\\\\
We start this section by formulating hypotheses related to the main results.
\begin{Hyp} \label{Hyp1}
The function $u$ satisfies $u\in C([0,T];W^{1,\infty}_0(\Omega))$.
\end{Hyp}
\begin{Hyp} \label{HypMM}
The nodal points of the finite element spaces $\Psi_h^0, \dots \Psi_h^{N_T}$ are given by a moving mesh.
\end{Hyp}
\begin{Hyp}\label{Hyp2}
The time increment $\Delta t$ satisfies the condition $\Delta t |u|_{C(W^{1,\infty})}\leq 1/8$.
\end{Hyp}

\begin{Rmk}
    Hypothesis \ref{Hyp2} is not a CFL condition since the mesh size $h$ is not included in the inequality. The time increment $\Delta t$ can be chosen independently of $h$.
\end{Rmk}

\begin{Hyp} \label{Hyp4}
The solution $\phi$ of problem (\ref{prob:c-d}) satisfies $\phi \in Z^3 \cap H^2 (0,T;H^2(\Omega)) \cap 
H^1 (0,T;H^3(\Omega))$
\end{Hyp}

%\\\\\\\\\\\\\\\\\\\\\\\\\\\\\\\\\\\\\\\\\\\\\\\\\\\\\\\\\\\\\\\\\\\
\subsection{Results regarding the Moving Mesh}\label{sec:resultsmm}
The first result we state concerns the moving mesh method~\eqref{MM_d}. 
\begin{Thm}[Non-overlapping condition for the moving mesh method] \label{Thm5}
Suppose that Hypothesis \ref{Hyp1} holds true. Let $C_0\in[0,1)$ be fixed, the set of nodal points $\{P^n_i:~i=1,\dots,N_p,~n=1,\dots,N_T\}$ be given by method~\eqref{MM_d}, and
\begin{equation}\label{eq:CFL}
    \Delta t |u|_{C^0(W^{1,\infty}(\Omega))}\leq C_0,
\end{equation}
then the set of nodal points describes a moving mesh, i.e., it holds that for any $n\in\{0,\dots,N_T\}$
\begin{equation} \label{th1}
    P^{n}_i < P^{n}_j; \quad i< j; \quad i,j\in\{1,\dots,N_p\}.
\end{equation}
\end{Thm}
\begin{Rmk}
Suppose that the nodal points of the finite element spaces $\Psi_h^0, \dots \Psi_h^{N_T}$ are governed by the moving mesh method~\eqref{MM_d} then Hypothesis~\ref{HypMM} is implied by condition \eqref{eq:CFL} due to Theorem~\ref{Thm5}.
\end{Rmk}
Next, we state two results that are necessary in order to derive the error estimates for the LGMM schemes. For $f\in C^0(\bar{\Omega}), t\in[0,T]$, and the time dependent P1-basis functions $\varphi_i(x, t)$ for $i\in\{1,\dots, N_p\}$ we define the time dependent Lagrange interpolation of $f$ by 
\begin{equation}
    \left[\Pi_h (t) f\right](x)\defeq\sum_{i=1}^{N_p} f(P_i(t))\psi_i(x,t).
\end{equation}
We also denote the difference operator $\bar{D}_{\Delta t} f \defeq \frac{f^n - f^{n-1}}{\Delta t}$.
\begin{Thm} \label{Prop1}
Let $\{\phi(t)=\phi(\cdot,t)\in\Psi ; t\in(0,T)\}$ be the solution of problem \eqref{prob:c-d}. Suppose that Hypothesis~\ref{HypMM} and Hypothesis~\ref{Hyp4} holds true. 
Then assuming $w \in C^0(W^{1,\infty}_0(\Omega))$ the following results hold.
\begin{enumerate}
\item[i)] There exists a positive constant $C=C(\|w \|_{C^0(L^\infty)})$ independent of $\Delta t$ and $h$ such that
  \begin{subequations}
\begin{equation} \label{eq:prop1a}
\left\Vert \bar{D}_{\Delta t}(\Pi^n_h\phi^n) - \frac{1}{\Delta t} \int_{t^{n-1}}^{t^{n}} \Pi_h(t) \frac{\partial\phi}{\partial t} (\cdot, t) dt \right \Vert 
\leq \frac{C h}{\sqrt{\Delta t}} \| \phi \| _{L^2(t^{n-1},t^n ; H^2(\Omega))}.
\end{equation}
\item[ii)] For a positive constant $C^\prime = C^\prime(\|w \|_{C^0(W^{1,\infty})})$ independent of $\Delta t$ and $h$ it holds
  \begin{equation} \label{eq:prop1b}
\left\Vert \bar{D}_{\Delta t}(\Pi^n_h\phi^n) - \frac{1}{\Delta t} \int_{t^{n-1}}^{t^{n}} \Pi_h(t) \frac{\partial\phi}{\partial t} (\cdot, t) dt \right \Vert _{\Psi^\prime}
\leq C h^2  \| \phi \| _{H^1(H^3)}.
\end{equation}
\end{subequations}
\end{enumerate}
\end{Thm}

\begin{Rmk}
    In case of a fixed mesh (velocity $w=0$), the corresponding $O(h^2)$-estimate of \eqref{eq:prop1b} can be obtained in $L^2$-norm.
\end{Rmk}

\begin{Cor} \label{Thm6}
Let $\{\phi(t)=\phi(\cdot,t)\in\Psi ; t\in(0,T)\}$ be the solution of problem \eqref{prob:c-d}. Suppose that Hypothesis~\ref{HypMM} and Hypothesis~\ref{Hyp4} hold true. Define $\eta(t)\defeq\phi (t)-\Pi_h (t) \phi(t)$. 
Then for $w \in C^0(W^{1,\infty}_0(\Omega))$ there exist positive constants $C=C(\|w \|_{C^0(L^\infty)})$ and $C^\prime=C^\prime(\|w \|_{C^0(W^{1,\infty})})$ independent of $\Delta t$ and $h$ such that the following bounds hold
\begin{subequations}
\begin{align} 
  \left\Vert \bar{D}_{\Delta t}\eta^n \right\Vert &\leq \frac{C h}{\sqrt{\Delta t}}  \| \phi \| _{H^1(t^{n-1},t^n;H^2(\Omega))}, \label{th2a}\\
  \left\Vert \bar{D}_{\Delta t}\eta^n \right\Vert_{\Psi^\prime}  
  &\leq C' h^2  \left( \frac{1}{\sqrt{\Delta t}} \| \phi \| _{H^1(t^{n-1},t^n;H^2(\Omega))} +  \| \phi \| _{H^1(H^3)} \right). 
  \label{th2b}
\end{align}
\end{subequations}
\end{Cor}

%\\\\\\\\\\\\\\\\\\\\\\\\\\\\\\\\\\\\\\\\\\\\\\\\\\\\\\\\\\\\\\\\\\\
\subsection{Results for the First-Order LGMM Scheme} 
In this section we state results for the scheme introduced in Section~\ref{sec:1storderLG}.

\begin{Cor}[Mass preserving property of the first-order LGMM Scheme] \label{Thm1}
  Suppose that Hypotheses \ref{Hyp1}, \ref{HypMM} and \ref{Hyp2} hold true. Let $\{\phi^n_h\}_{n=1}^{N_T}$ be the solution of the numerical scheme \eqref{LG} for a given initial datum $\phi_h^0$. Then it holds for $n=0,1,\dots,N_T$ that
\begin{equation} \label{th3}
\int_{\Omega}\phi^n_hdx = \int_{\Omega}\phi^0_h dx + \Delta t \sum_{i=1}^{n}\left(\int_{\Omega}f^i dx + \int_{\Gamma} g^i ds \right).
\end{equation}
\end{Cor}

\begin{Cor}[Stability of the first-order LGMM scheme] \label{Thm2}
Suppose that Hypotheses \ref{Hyp1}, \ref{HypMM} and \ref{Hyp2} hold true. Let $F\in H^1(0,T;\Psi')$ be given. For the given function $\phi^0_h \in \Psi_h$, let $\{\phi^n_h\}^{N_T}_{n=1} \subset \Psi_h$ be the numerical solutions of scheme \eqref{LG}. Then there exists a constant $C>0 $ independent of $h$ and $\Delta t$ such that 
\begin{equation} \label{th4}
 \| \phi_h \| _{\ell^\infty(L^2)} + \sqrt{\nu} \| \nabla \phi_h\| _{\ell^2(L^2)} \leq  C \left( \| \phi^0_h \| + \| F \| _{\ell^2(\Psi'_h)}\right).
\end{equation}
\end{Cor}

\begin{Rmk}
  The mass-preserving and stability properties of the first-order Lagrange-Galerkin scheme with fixed mesh (Theorem 1 and Theorem 2 of \cite{rui2010massconsercharac}) are maintained in the first-order LGMM scheme.
\end{Rmk}

\begin{Thm}[Error estimates for the first-order LGMM scheme] \label{Thm7}
  Suppose that Hypotheses \ref{Hyp1}, \ref{HypMM}, \ref{Hyp2}, and \ref{Hyp4} hold true. Let $F\in H^1(0,T;\Psi')$ be given. Assuming the initial datum $\phi^0_h=\Pi^0_h\phi^0\in \Psi_h$ let $\{\phi(t)=\phi(\cdot,t)\in\Psi ; t\in(0,T)\}$ be the solution of problem \eqref{prob:c-d} and $\{\phi^n_h\}_{n=1}^{N_T}$ the numerical solutions of scheme \eqref{LG}. Then there exists a constant $C>0 $ independent of $h$ and $\Delta t$ such that 
\begin{equation} \label{th5}
\begin{split}
 \| \phi_h - \phi \| _{\ell^\infty(L^2)} + \sqrt{\nu}   \| \nabla (\phi_h-\phi)\| _{\ell^2(L^2)} &\leq   C  (\Delta t + h^2)  \| \phi \| _{Z^2 \cap H^1(H^2) \cap H^1(H^3)}.
\end{split}
\end{equation}
\end{Thm}
\begin{Rmk}\label{rmk:lowerorder} Using the bound \eqref{th2a} instead of \eqref{th2b} in the proof of Theorem~\ref{Thm7} a first order bound that requires lower regularity of $\phi$ is obtained. Namely, under the assumptions of Theorem~\ref{Thm7} there exists a constant $C>0$ independent of $h$ and $\Delta t$ such that 
\begin{equation*} 
 \| \phi_h - \phi \| _{\ell^\infty(L^2)} + \sqrt{\nu}   \| \nabla (\phi_h-\phi)\| _{\ell^2(L^2)} \leq   C  (\Delta t + h)  \| \phi \| _{Z^2 \cap H^1(H^2) \cap  H^1(H^2)}.
\end{equation*}
\end{Rmk}
% \\\\\\\\\\\\\\\\\\\\\\\\\\\\\\\\\\\\\\\\\\\\\\\\\\\\\\\\\\\\\\\\\\\
\subsection{Results for the Second-Order LGMM Scheme}
The results in this section concern the second-order scheme introduced in Section~\ref{sec:2ndorderLG}.
\begin{Cor}[Mass preserving property of the second-order LGMM scheme] \label{Thm3}
Suppose that Hypotheses \ref{Hyp1}, \ref{HypMM} and \ref{Hyp2} hold true. Let $\{\phi^n_h\}_{n=1}^{N_T}$ be the solution of the numerical scheme~\eqref{LG2} for a given initial datum $\phi_h^0$. Then it holds for $n=1,2,\dots,N_T$ that 
\begin{equation} \label{th6}
\int_{\Omega}\left( \frac{3}{2} \phi^n_h - \frac{1}{2} \phi^{n-1}_h \right) dx = \int_{\Omega}\phi^0_h dx + \Delta t \sum_{i=1}^{n}\left(\int_{\Omega}f^i dx + \int_{\Gamma} g^i ds \right).
\end{equation}
\end{Cor}

\begin{Cor}[Stability for the Second-Order LGMM Scheme] \label{Thm4}
 Suppose that Hypotheses \ref{Hyp1}, \ref{HypMM} and \ref{Hyp2} hold true. Let $F\in H^1(0,T;\Psi')$ be given. For a given function $\phi^0_h \in \Psi_h$ let $\{\phi^n_h\}^{N_T}_{n=1} \subset \Psi_h$ be the numerical solutions of scheme~\eqref{LG2}. Then, there exists a constant $C>0$ independent of $h$ and $\Delta t$ such that 
\begin{equation} \label{th7}
     \| \phi_h \| _{\ell^\infty(L^2)} + \sqrt{\nu} \| \nabla \phi_h\| _{\ell^2(L^2)} \leq C \left( \| \phi^0_h \|   +  \| F \| _{\ell^2(\Psi'_h)}\right).
\end{equation}
\end{Cor}

\begin{Rmk}
Also in case of the second-order Lagrange--Galerkin scheme the mass-preserving and stability properties of the fixed mesh method (Theorem 1 and Theorem 2 of \cite{futai2022masspresertwo}) are maintained in the LGMM scheme.
\end{Rmk}

\begin{Thm}[Error Estimates the Second-Order LGMM Scheme] \label{Thm8}
Suppose that Hypotheses \ref{Hyp1}, \ref{HypMM}, \ref{Hyp2}, and \ref{Hyp4} hold true. Let $F\in H^1(0,T;\Psi')$ be given. For a given function $\phi^0_h=\Pi^0_h\phi^0\in \Psi_h$ let $\{\phi(t)=\phi(\cdot,t)\in\Psi ; t\in(0,T)\}$ be the solution of problem \eqref{prob:c-d} and $\{\phi^n_h\}_{n=1}^{N_T}$ be the numerical solutions of scheme (\ref{LG2}). Then, there exists a constant $C>0 $ independent of $h$ and $\Delta t$ such that 
\begin{equation} \label{th8}
\begin{split}
 \| \phi_h - \phi \| _{\ell^\infty(L^2)} + \sqrt{\nu}   \| \nabla (\phi_h-\phi)\| _{\ell^2(L^2)} &\leq  C (\Delta t^2 + h^2)  \| \phi \| _{Z^3 \cap H^2(H^2)\cap  H^1(H^3)}.
\end{split}
\end{equation}
 \end{Thm}
\begin{Rmk} In analogy to Remark~\ref{rmk:lowerorder} in the proof of Theorem~\ref{Thm8} the bound \eqref{th2a} can be used instead of \eqref{th2b} to obtain a first order bound that requires lower regularity of $\phi$. Namely, under the assumptions of Theorem~\ref{Thm8} there exists a constant $C>0$ independent of $h$ and $\Delta t$ such that 
\begin{equation*} 
 \| \phi_h - \phi \| _{\ell^\infty(L^2)} + \sqrt{\nu}   \| \nabla (\phi_h-\phi)\|_{\ell^2(L^2)} \leq   C  (\Delta t^2 + h) \| \phi \|_{Z^3 \cap H^2(H^2)}
\end{equation*}
\end{Rmk}
% \\\\\\\\\\\\\\\\\\\\\\\\\\\\\\\\\\\\\\\\\\\\\\\\\\\\\\\\\\\\\\\\\\\
 \section{PROOFS}
 In this section we provide proofs for the results stated in Section~\ref{sec:results}.
%\\\\\\\\\\\\\\\\\\\\\\\\\\\\\\\\\\\\\\\\\\\\\\\\\\\\\\\\\\\\\\\\\\\
\subsection{Proofs of the Results regarding the Moving Mesh}
\subsubsection{Proof of Theorem \ref{Thm5}}
We show property~\eqref{th1} inductively. Hence, suppose $P^{n-1}_i < P^{n-1}_j;$ $ i < j;$ $ i,j\in\{1,\dots,N_p\}$, we show that~\eqref{th1} holds true.
Let $h^{n-1}_i\defeq P^{n-1}_{i+1}-P^{n-1}_i$ for $i \in \{1,\dots,N_p-1\}$. It is sufficient to show that $h^{n}_i > 0$ for $i \in \{1,\dots,N_p-1\}$. Shifting the index $i$ in scheme \eqref{mm_d_1}, we have
\begin{equation} \label{mm_d_4}
 \frac{P^{n}_{i+1}-P^{n-1}_{i+1}}{\Delta t}=u^{n-1}(P^{n-1}_{i+1})+\nu\frac{P^{n}_{i+2}-2P^{n}_{i+1}+P^{n}_{i}}{(P^{n-1}_{i+1}-P^{n-1}_{i})(P^{n-1}_{i+2}-P^{n-1}_{i+1})}, \quad i= 1,\dots,N_p-2.
\end{equation}
By subtracting~\eqref{mm_d_1} from~\eqref{mm_d_4} we obtain
\begin{equation*}
    \frac{h^{n}_{i}-h^{n-1}_{i}}{\Delta t}=u^{n-1}(P^{n-1}_{i+1})-u^{n-1}(P^{n-1}_{i})+\nu\left[\frac{h^{n}_{i+1}-h^{n}_{i}}{h^{n-1}_{i}h^{n-1}_{i+1}}-\frac{h^{n}_{i}-h^{n}_{i-1}}{h^{n-1}_{i-1}h^{n-1}_{i}}\right], \quad i= 2,\dots,N_p-2.
\end{equation*}
Rearranging the terms it follows for $i=2,\dots,N_p-2$ 
\begin{align}\label{eq:hsystem}
    \left[\frac{1}{\Delta t}+\nu\left(\frac{1}{h^{n-1}_{i}h^{n-1}_{i+1}}+\frac{1}{h^{n-1}_{i-1}h^{n-1}_{i}}\right)\right]h^{n}_i&-\nu\frac{1}{h^{n-1}_{i}h^{n-1}_{i+1}}h^{n}_{i+1}-\nu\frac{1}{h^{n-1}_{i-1}h^{n-1}_{i}}h^{n}_{i-1} \nonumber\\ 
                                                                                                                                &=u^{n-1}(P^{n-1}_{i+1})-u^{n-1}(P^{n-1}_{i})+\frac{h^{n-1}_i}{\Delta t}.
\end{align}
Using \eqref{eq:CFL} we derive a lower bound of the right hand side in~\eqref{eq:hsystem} as follows:
\begin{align}
  u^{n-1}(P^{n-1}_{i+1})-u^{n-1}(P^{n-1}_{i})+\frac{h^{n-1}_i}{\Delta t}   &\geq \frac{h^{n-1}_i}{\Delta t}-|u^{n-1}(P^{n-1}_{i+1})-u^{n-1}(P^{n-1}_{i})| \nonumber \\
    &\geq \frac{h^{n-1}_i}{\Delta t}-h^{n-1}_i|u|_{C^0(W^{1,\infty}_0)}=\frac{h^{n-1}_i}{\Delta t}\left(1-\Delta t|u|_{C^0(W^{1,\infty}_0)}\right) \nonumber \\
    &>0. \label{th1.1}
\end{align}
Note that from~\eqref{mm_d_1} for $i=2$ we have 
\begin{equation*}
    \frac{P^{n}_2-P^{n-1}_2}{\Delta t}=u^{n-1}(P^{n-1}_2)+\nu\frac{h^{n}_2-h^{n}_1}{h^{n-1}_1h^{n-1}_2},
\end{equation*}
and from~\eqref{mm_d_2} follows $P^{n-1}_1=a$, which implies $ \frac{P^{n}_1-P^{n-1}_1}{\Delta t}=0$. 
Therefore, it holds
\begin{equation}\label{eq:hsystemleft}
    \left(\frac{1}{\Delta t}+\nu\frac{1}{h^{n-1}_{1}h^{n-1}_{2}}\right)h^{n}_1-\nu\frac{1}{h^{n-1}_{1}h^{n-1}_{2}}h^{n}_{2}    =u^{n-1}(P^{n-1}_{2})-u^{n-1}(P^{n-1}_{1})+\frac{h^{n-1}_1}{\Delta t}.
\end{equation}
Similarly, from~\eqref{mm_d_1} for $i=N_p-1$ we have 
\begin{equation*}
    \frac{P^{n}_{N_p-1}-P^{n-1}_{N_p-1}}{\Delta t}=u^{n-1}(P^{n-1}_{N_p-1})+\nu\frac{h^{n}_{N_p-1}-h^{n}_{N_p-2}}{h^{n-1}_{N_p-2}h^{n-1}_{N_p-1}}
\end{equation*}
and from \eqref{mm_d_2} we obtain $P^{n-1}_{N_p}=b$, which implies $ \frac{P^{n}_{N_p}-P^{n-1}_{N_p}}{\Delta t}=0$. 
Therefore, we have
\begin{equation}\label{eq:hsystemright}
  \left(\frac{1}{\Delta t}+\nu\frac{1}{h^{n-1}_{N_p-2}h^{n-1}_{N_p-1}}\right)h^{n}_{N_p-1} -\nu\frac{1}{h^{n-1}_{N_p-2}h^{n-1}_{N_p-1}}h^{n}_{N_p-2} =u^{n-1}(P^{n-1}_{N_p})-u^{n-1}(P^{n-1}_{N_p-1})+\frac{h^{n-1}_{N_p-1}}{\Delta t}.
\end{equation}
Proceeding as in \eqref{th1.1} positivity of the right hand sides in both \eqref{eq:hsystemleft} and \eqref{eq:hsystemright} follows. Combining \eqref{eq:hsystemleft}, \eqref{eq:hsystem} and \eqref{eq:hsystemright} yields a linear system with unknown variables $h^{n}_1$, \dots, $h^n_{N_p-1}$ and a strictly diagonally dominant coefficient matrix, which is thus an M-matrix. Since an M-matrix $A$ has the property that $Ax>0$ implies $x>0$, the solution of the linear system is positive, i.e.,  $h^{n}_{1}$, \dots,$h^n_{N_p-1} > 0$, hence~\eqref{th1} follows. \qed
 
%//////////////////////////////////////////////////////////////////
\subsubsection{Proof of Theorem \ref{Prop1}}
First, we state the following lemma, which plays an important role in the proof of Theorem \ref{Prop1}. The proof of the lemma is given in Appendix~\ref{sec:proofl1}.
\begin{Lem}\label{Lem3}
Let $\{\phi(t)=\phi(\cdot,t)\in\Psi ; t\in(0,T)\}$ be the solution of problem~\eqref{prob:c-d} and suppose that Hypothesis~\ref{HypMM} holds true. For $x\in\bar{\Omega}$ and $t\in[t^{n-1},t^n]$ we define
\begin{equation*}
I (x,t) \defeq \sum_{i=1}^{N_p}\phi(P_i(t),t)\left[\frac{\partial}{\partial t}\psi_i(x,t)\right],
\end{equation*}
where $P_i(t)$ for $i\in\{1,\dots,N_p\}$ are the nodal point positions defined in \eqref{eq:movingpoints} and $\psi_i(x,t)$ for $i\in\{1,\dots,N_p\}$ denote the time extended P1 basis functions. We assume $x\in [P_k(t), P_{k+1}(t)]$ for a $k\in{\{1,\dots,N_p -1\}}$. Then $I (x,t)$ can be expressed as 
\begin{equation} \label{claim1}
I(x,t)=-\frac{\phi(P_{k+1}(t),t)-\phi(P_k(t),t)}{P_{k+1}(t)-P_k(t)}\left[w^n(P_{k+1}(t))\psi_{k+1}(x,t)+w^n(P_k(t))\psi_k(x,t))\right].
\end{equation}
\end{Lem}

\paragraph{Proof of Theorem \ref{Prop1}} We first define the interpolation operators
\begin{align*}
\left(\Pi^\ell_h\phi^\ell\right) (x) &\defeq \sum_{i=1}^{N_p} \phi^\ell (P_i^\ell)\psi^\ell_i (x), \qquad \ell \in \{ n-1, n\}.
\end{align*}
Then we rewrite their difference as 
\begin{align}
&\left(\Pi^n_h\phi^n-\Pi^{n-1}_h\phi^{n-1}\right)(x) \nonumber \\
&=\int_{t^{n-1}}^{t^n}\frac{d}{dt} (\Pi_h(t)\phi(\cdot,t))(x) dt \nonumber \\
&=\sum_{i=1}^{N_p}\int_{t^{n-1}}^{t^n}\frac{\partial}{\partial t} [\phi(P_i(t),t)\psi_i(x,t)] dt \nonumber \\
&=\sum_{i=1}^{N_p}\int_{t^{n-1}}^{t^n}\left(\left[\frac{\partial}{\partial t}\phi(P_i(t),t)\right]\psi_i(x,t)+\phi(P_i(t),t)\left[\frac{\partial}{\partial t}\psi_i(x,t)\right]\right) dt\nonumber \\
&=\sum_{i=1}^{N_p}\int_{t^{n-1}}^{t^n}\left(\left[\frac{\partial\phi}{\partial t}(P_i(t),t) + \frac{d P_i}{d t}(t)(\nabla\phi)(P_i(t),t)\right]
\psi_i(x,t)+\phi(P_i(t),t)\left[\frac{\partial}{\partial t}\psi_i(x,t)\right]\right) dt\nonumber \\
&=\sum_{i=1}^{N_p}\int_{t^{n-1}}^{t^n}\left(\left[\frac{\partial\phi}{\partial t}(P_i(t),t) + w^n(P_i(t))(\nabla\phi)(P_i(t),t)\right]
\psi_i(x,t)+\phi(P_i(t),t)\left[\frac{\partial}{\partial t}\psi_i(x,t)\right]\right) dt\nonumber \\
&=\int_{t^{n-1}}^{t^n} \Pi_h(t) \left[\frac{\partial\phi}{\partial t} (\cdot,t) + w^n(\cdot) \nabla\phi(\cdot,t)\right] (x) dt + \int_{t^{n-1}}^{t^n} I(x,t) dt. \label{R2}
\end{align}

The rest of the proof concerns the last integral in \eqref{R2}. Without loss of generality let $t\in[t^{n-1},t^n]$ and $x\in K_k(t) \defeq [P_k(t), P_{k+1}(t)]$. For brevity we introduce the notations:
\begin{align*}
w^n_{k}&\defeq w^n(P^{n-1}_{k}), \qquad
w^n_{k+1}\defeq w^n(P^{n-1}_{k+1}), \qquad
  h_k = P_{k+1}(t) - P_k(t), \\
  \phi_{k}&\defeq \phi(P_{k}(t),t), \qquad
  \phi_{k+1}\defeq \phi(P_{k+1}(t),t).
\end{align*}
We're now in the position to show i). 
Due to the Taylor expansions
 \begin{align*}
\phi_{k+1} - \phi_k
&=\phi(P_{k}(t)+h_k(t),t) - \phi(P_{k}(t),t)\\
&=h_k(t) ( \nabla\phi ) (P_k(t),t) + h_k^2(t) \int_0^1 \int_0^{s_0} (\nabla^2 \phi)(P_k(t)+s_1 h_k(t), t) ds_1 ds_0, \\
\phi_{k+1} - \phi_k 
&=\phi(P_{k+1}(t),t) - \phi(P_{k+1}(t)-h_k(t),t)\\
&=h_k(t) ( \nabla\phi ) (P_{k+1}(t),t) - h_k^2(t) \int_0^1 \int_0^{s_0} (\nabla^2 \phi)(P_{k+1}(t)-s_1 h_k(t), t) ds_1 ds_0. 
 \end{align*}
we obtain the identity
 \begin{align}\label{eq:taylorphidiff}
&\frac{\phi_{k+1} - \phi_k}{P_{k+1}(t)-P_k(t)}\left[ w^n_{k}\psi_{k}(x,t)+w^n_{k+1}\psi_{k+1}(x,t)\right] \notag \\
&\quad= \left[ \Pi_h(t) w^n (\cdot) \nabla\phi (\cdot,t) \right] (x) +  h_k(t)  w^n_{k}\psi_{k}(x,t) \int_0^1 \int_0^{s_0}  (\nabla^2 \phi)(P_k(t)+s_1 h_k(t), t)  ds_1 ds_0   \notag \\
&\quad - h_k(t) w^n_{k+1}\psi_{k+1}(x,t) \int_0^1 \int_0^{s_0} (\nabla^2 \phi)(P_{k+1}(t) - s_1 h_k(t), t) \, ds_1 ds_0 .
 \end{align}
By using Lemma~\ref{Lem3}, we substitute~\eqref{eq:taylorphidiff} into~\eqref{R2}, and through a change of variable, we proceed to compute
\begin{align*}
\left\vert \left( \Pi^n_h\phi^n - \Pi^{n-1}_h\phi^{n-1} \right) - \int_{t^{n-1}}^{t^n} \left( \Pi_h(t) \frac{\partial\phi}{\partial t} (\cdot,t) \right) dt \right\vert  &= \left\vert \int_{t^{n-1}}^{t^n}\Pi_h(t) w^n (\cdot) \nabla\phi (\cdot,t) (x) + \int_{t^{n-1}}^{t^n}I(x,t) dt \right\vert \\
&\leq c \|w \|_{C^0(L^\infty)} \int_{t^{n-1}}^{t^n} \int_{P_k(t)}^{P_{k+1}(t)} |(\nabla^2 \phi)(x , t)| dx dt.
\end{align*}
By taking the $L^2$-norm over $\Omega$ and applying the Cauchy--Schwartz inequality on the right hand side, we obtain
\begin{align}\label{eq:prop1ifinal}
    \left\Vert \left( \Pi^n_h\phi^n - \Pi^{n-1}_h\phi^{n-1} \right) - \int_{t^{n-1}}^{t^n} \left( \Pi_h(t) \frac{\partial\phi}{\partial t} (\cdot,t) \right) dt\right \Vert^2 
    &\leq c^2 h \Delta t \|w \|^2_{C^0(L^\infty)}  \int_{\Omega} \int_{t^{n-1}}^{t^n} \int_{K_{\tilde k}(t)} (\nabla^2 \phi)(y , t)^2 dy dt  dx \notag \\
    &= c^2 h \Delta t \|w \|^2_{C^0(L^\infty)} \int_{t^{n-1}}^{t^n} \sum_{k=1}^{N_p-1} h_{k(t)} \int_{K_k(t)} (\nabla^2 \phi)(y , t)^2 dy  dt \notag \\
    &\leq c^2 h^2 \Delta t \|w \|^2_{C^0(L^\infty)} ||\nabla^2\phi||^2_{L^2(t^{n-1},t^n;L^2(\Omega))},
\end{align}
where the dynamic index $\tilde k=\tilde k(x,t)$ is defined such that $x\in K_{\tilde k}(t)$. 
To complete the proof, we take the square root and divide both sides of \eqref{eq:prop1ifinal} by $\Delta t$, obtaining \eqref{eq:prop1a}.

Next, we proof ii). Therefore, we first rewrite and then further expand the last double integral in~\eqref{eq:taylorphidiff} as follows
\begin{align*}
\int_0^1 \int_0^{s_0} (\nabla^2 \phi)(P_k(t) + (1-s_1) h_k(t), t) ds_1 ds_0 &= \int_0^1 \int_0^{s_0}  (\nabla^2 \phi)(P_k(t)+s_1 h_k(t), t)  ds_1 ds_0  \\
&\quad + h_k(t) \int_0^1 \int_0^{s_0} \int_{s_1}^{1-s_1}  (\nabla^3 \phi)(P_k(t)+s_2 h_k(t), t)  ds_2  ds_1 ds_0.
\end{align*}
Additionally we introduce the Taylor expansion
 \begin{align*}
w^n_{k+1} &= w^n_{k} + h_k \int_0^1 (\nabla w^n)(P_k(t) + sh_k(t)) ds.
 \end{align*}
By using Lemma~\ref{Lem3} and substituting the above expressions into \eqref{R2} we compute 
\begin{align*}
    A^n(x)&:= \frac{1}{\Delta t} \left[ \left( \Pi^n_h\phi^n-\Pi^{n-1}_h\phi^{n-1} \right) (x) - \int_{t^{n-1}}^{t^n} \left( \Pi_h(t) \frac{\partial\phi}{\partial t} (\cdot,t) \right) (x) dt \right] \\ 
    &= -\frac{1}{\Delta t} \int_{t^{n-1}}^{t^n} h_k(t) w^n_k(\psi_k(x,t) - \psi_{k+1}(x,t)) \int_0^1  \int_0^{s_0}  (\nabla^2 \phi)(P_k(t)+s_1 h_k(t), t)  ds_1 ds_0  dt \\
    &\quad + \frac{1}{\Delta t} \int_{t^{n-1}}^{t^n}  h^2_k(t) \psi_{k+1}(x,t) w^n_k \int_0^1 \int_0^{s_0} \int_{s_1}^{1-s_1}  (\nabla^3 \phi)(P_k(t)+s_2 h_k(t), t)  ds_2 ds_1 ds_0  dt \\
          &\quad + \frac{1}{\Delta t} \int_{t^{n-1}}^{t^n} \left( h^2_k(t) \psi_{k+1}(x,t) \int_0^1 (\nabla w^n)(P_k(t) + sh_k(t)) ds \right. \\
  & \quad \qquad \qquad \left. \int_0^1 \int_0^{s_0} (\nabla^2 \phi)(P_k(t) + (1-s_1) h_k(t), t) ds_1 ds_0\right) dt \\ 
    &\eqdef A^n_1(x) + A^n_2(x) + A^n_3(x).
\end{align*}
We proceed by estimating the ($\Psi^\prime$)-norm, i.e., $\|A^n\|_{\Psi^\prime}\leq \|A^n_1\|_{\Psi^\prime} + \|A^n_2\|_{\Psi^\prime} + \|A^n_3\|_{\Psi^\prime}$. The following bounds hold
\begin{align}
   \|A^n_1\|_{\Psi^\prime} &\leq c_1 h^2 \|w\|_{C^0(L^\infty)} \|\phi\|_{H^1(H^3)} \label{eq:a1bound}\\
   \|A^n_2\|_{\Psi^\prime} &\leq c_2 h^2 \|w\|_{C^0(L^\infty)} \|\phi\|_{H^1(H^3)} \label{eq:a2bound}\\
   \|A^n_3\|_{\Psi^\prime} &\leq c_3 h^2 \|w\|_{C^0(W^{1,\infty})} \|\phi\|_{H^1(H^3)} \label{eq:a3bound}
\end{align}
and their detailed proofs are provided in Appendix~\ref{sec:cor2}. Combining all the bounds, we obtain the estimate
\begin{align*}
   \|A\|_{\Psi^\prime} &\leq C h^2 \|w\|_{C^0(W^{1,\infty})} \|\phi\|_{H^1(H^3)},
\end{align*}
which completes the proof of \eqref{eq:prop1b}. \qed

\subsubsection{Proof of Corollary \ref{Thm6}}
We first recall an error estimate for the Lagrange interpolation that follows from \cite[Theorem 4.4.20]{brenner2008}.
\begin{Lem}\label{lem:interpolation} We suppose that Hypothesis~\ref{HypMM} holds true and fix $t\in [0, T]$. Let $\Pi_h=\Pi_h(t)= \sum_{i=1}^{N_p} f(P_i(t))\psi_i(x,t)$ be the Lagrange interpolation operator at time $t$ and $v \in H^2(\Omega)$. Then there exists a constant $C>0$ independent of $h$ such that
  \begin{equation}\label{eq:interpolationerror}
    \| \Pi_h v - v \|_{H^s(\Omega)} \leq C h^{2-s} | v |_{H^2(\Omega)} \qquad \text{for }s \in \{0, 1\}.
  \end{equation}
\end{Lem}
\paragraph{Proof of Corollary \ref{Thm6}} We show \eqref{th2a} and assume $t\in[t^{n-1},t^n]$. 
By applying the bound \eqref{eq:prop1a} from Theorem \ref{Prop1}, Cauchy-Schwartz's inequality, as well as Lemma~\ref{lem:interpolation}, we obtain the estimate, 
\begin{align}
\left\Vert\eta^n - \eta^{n-1} \right\Vert 
&=\left\Vert(\phi^n-\Pi^{n}_h\phi^n)-(\phi^{n-1}-\Pi^{n-1}_h\phi^{n-1}) \right\Vert  \nonumber \\
&=\left\Vert\phi^n-\phi^{n-1} - \Pi^n_h\phi^n-\Pi^{n-1}_h\phi^{n-1} \right\Vert \nonumber\\ 
&\leq \left\Vert \int_{t^{n-1}}^{t^n}\left(\frac{\partial\phi}{\partial t} - \Pi_h(t) \frac{\partial\phi}{\partial t} \right) dt \right\Vert + c_1 h \sqrt{\Delta t}  \| \phi \| _{L^2(t^{n-1},t^n);H^2(\Omega))} \nonumber\\
&\leq \sqrt{\Delta t} \sqrt{\int_{t^{n-1}}^{t^n}\left\Vert\frac{\partial\phi}{\partial t} - \Pi_h \frac{\partial\phi}{\partial t}\right\Vert^2_{L^2(\Omega)} dt } +  c_1 h \sqrt{\Delta t}  \| \phi \| _{L^2(t^{n-1},t^n);H^2(\Omega))}  \nonumber\\
& \leq c_2  \sqrt{\Delta t} \left( h^2 \sqrt{\int_{t^{n-1}}^{t^n}\left|\frac{\partial\phi}{\partial t}\right|^2_{H^2(\Omega)} dt } + h  \| \phi \| _{L^2(t^{n-1},t^n);H^2(\Omega))} \right) \nonumber \\
&\leq c_2 h  \sqrt{\Delta t} \left( h   \| \phi \| _{H^1(t^{n-1},t^n;H^2(\Omega))} +  \| \phi \| _{L^2(t^{n-1},t^n);H^2(\Omega))} \right). \label{Deltaeta}
\end{align}
To complete the proof, we divide both sides of \eqref{Deltaeta} by $\Delta t$ and obtain \eqref{th2a}. The bound \eqref{th2b} is obtained repeating the above estimates in the $\Psi^\prime$-norm, using \eqref{eq:prop1b} instead of \eqref{eq:prop1a} and embedding $L^2$ in $\Psi^\prime$. \qed
 
%\\\\\\\\\\\\\\\\\\\\\\\\\\\\\\\\\\\\\\\\\\\\\\\\\\\\\\\\\\\\\\\\\
\subsection{Proofs of the Results for the First-Order LGMM Scheme}\label{sec:proofslg1}
\paragraph{Proofs of Corollary \ref{Thm1} and Corollary \ref{Thm2}} The proof of the two corollaries follows directly from Theorem 1 and Theorem 2 in \cite{rui2010massconsercharac}, respectively. For convenience we provide the proofs in Appendix~\ref{sec:proofmasslg1} and Appendix~\ref{sec:proofstabilitylg1}.

To prove Theorem \ref{Thm7}, we first state the following lemma.
\begin{Lem}[Evaluation of composite functions \cite{futai2022masspresertwo,rui2010massconsercharac}] \label{lem:composite}
Let $a$ be a function in $W^{1,\infty}_0(\Omega)^d$ satisfying $\Delta t  \| a \| _{1,\infty}\leq 1/4$ and consider the mapping $X_1(a,\Delta t)$ defined in (\ref{LG}). Then, the following inequalities hold.
\begin{subequations}
    \begin{align}
     \| \psi\circ X_1(a,\Delta t) \|  &\leq (1+c_1\Delta t)  \| \psi \| , &&\forall \psi \in L^2(\Omega), \label{eq:compboundL2nodiff}\\
     \| \psi - \psi\circ X_1(a,\Delta t) \|  &\leq c_2\Delta t  \| \psi \| _{H^1(\Omega)}, &&\forall \psi \in H^1(\Omega), \label{eq:compboundL2}\\
     \| \psi - \psi\circ X_1(a,\Delta t) \| _{H^{-1}(\Omega)} &\leq c_3\Delta t  \| \psi \| , &&\forall \psi \in L^2(\Omega). \label{eq:compboundH-1}
\end{align}
\end{subequations}
\end{Lem}

\paragraph{Proof of Theorem \ref{Thm7}} We define the terms
\begin{align*}
e^n_h \defeq \phi ^n_h - \Pi^n_h \phi^n, \qquad
\eta(t) \defeq \phi (t) - \Pi_h (t) \phi(t).
\end{align*}
By substituting the error $e^n_h$ in the numerical scheme (\ref{LG}), we obtain the following expression:
\begin{equation} \label{err}
\left(\frac{e^n_h-[e^{n-1}_h\circ X^n_1]\gamma^n}{\Delta t},\psi_h\right)+a_0(e^n_h,\psi_h)=\langle R_h^n,\psi_h\rangle,  \hspace{0.2cm} \forall\psi_h\in \Psi^n_h,
\end{equation}
where the residual on the right hand side is given by
\begin{align*}
R^n_h & \defeq R^n_1 + R^n_2 + R^n_3, \\
R^n_1 & \defeq \frac{\partial \phi^n}{ \partial t} + \nabla \cdot \ (u^n \phi^n) - \frac{\phi^n - [\phi^{n-1}\circ X^n_1]\gamma^n}{\Delta t}, \\
R^n_2 & \defeq \frac{\eta^n - [\eta^{n-1}\circ X^n_1]\gamma^n}{\Delta t}, \\
\langle R^n_3, \psi_h\rangle &\defeq a_0(\eta^n,\psi_h).
\end{align*}
To obtain an estimate on $\|R_1\|$, we follow the error estimate framework for the  convection-diffusion problem on a static mesh (details are given in Appendix~\ref{sec:boundr11}), which gives us
\begin{align}
 \| R_1 \| _{\ell^2(\Psi'_h)} &\leq  c_4 \Delta t  \| \phi \| _{Z^2(0,T)} \label{eq:r1bound}.
\end{align}
In case of linear elements in one dimension that are considered here we have $R_3^n=0$ as is shown in Appendix~\ref{sec:boundr3}.
To compute a bound for $R^n_2$ we rewrite it as 
\begin{align*}
  R^n_2 &= \frac{\eta^n - [\eta^{n-1}\circ X^n_1]\gamma^n}{\Delta t} \\
  &= \frac{\eta^n - \eta^{n-1}}{\Delta t} + \frac{\eta^{n-1} - \eta^{n-1}\circ X^n_1}{\Delta t} + \frac{(\eta^{n-1} \circ X^n_1 )(1-\gamma^n)}{\Delta t}.
\end{align*}
Then, using \eqref{eq:compboundH-1} and \eqref{eq:compboundL2nodiff}, noting that thanks to Hypothesis~\ref{Hyp1} it holds $1-\gamma^n \leq c_2 \Delta t$, employing \eqref{th2b}, Lemma~\ref{lem:interpolation} and embedding $L^2(\Omega)$ in $H^1(\Omega)^\prime$, we obtain the following
\begin{align}
 \| R^n_2 \|_{(\Psi^n_h)'}  &\leq \left\Vert\frac{\eta^n - \eta^{n-1}}{\Delta t}\right\Vert_{(\Psi^n_h)'} + c_3  \| \eta^{n-1} \| + c_6  \|  \eta^{n-1} \circ X^n_1  \|  \notag \\
&\leq c_7 \left[\frac{h^2}{\sqrt{\Delta t}}  \| \phi \| _{H^1(t^{n-1},t^n;H^2(\Omega))} + h^2 \| \phi \| _{H^1(H^3)} +   \|  \eta^{n-1} \|  \right] \notag \\
          &\leq c_8 \left[ \frac{h^2}{\sqrt{\Delta t}}  \| \phi \| _{H^1(t^{n-1},t^n;H^2(\Omega))} +  h^2 \| \phi \| _{H^1(H^3)} \right]  \label{eq:r2bound_0}
\end{align}
Hence, by taking the $\ell^2$-norm
\begin{align}
  \|R_2\|_{\ell^2(\Psi'_h)} 
&\leq c_9 h^2 \left(\| \phi \| _{H^1(0,T;H^2(\Omega))} +  \| \phi \| _{H^1(H^3)} \right). \label{eq:r2bound}
\end{align}
By combining the estimates \eqref{eq:r1bound}, \eqref{eq:r2bound}, and taking into account the fact $ R^n_3 =0$, we get
\begin{equation}\label{eq:rbound}
 \| R_h \| _{\ell^2(\Psi'_h)} \leq C \left( \Delta t  \| \phi \| _{Z^2(0,T)} + h^2  \| \phi \| _{H^1(0,T;H^2(\Omega))}  + h^2  \| \phi \|_{H^1(H^3)} \right) ,
\end{equation}
as estimate for the total residual, where $C>0$ is independent of $h$ and $\Delta t$. Applying the stability result from Theorem~\ref{Thm2} to \eqref{err} and using the bound \eqref{eq:rbound} we obtain the error estimate~\eqref{th5}. \qed
%\\\\\\\\\\\\\\\\\\\\\\\\\\\\\\\\\\\\\\\\\\\\\\\\\\\\\\\\\\\\\\\\\\\\\\
\subsection{Proofs of the Results for the Second-Order LGMM Scheme}\label{sec:proofslg2}
\paragraph{Proofs of Corollary \ref{Thm3} and Corollary \ref{Thm4}} The proof of both corollaries follows from Theorem 1 and Theorem 2 in \cite{futai2022masspresertwo}, respectively.  For convenience we provide the proofs in Appendices~\ref{sec:proofmasslg2} and \ref{sec:proofstabilitylg2}. \\

Next, we state the following lemma which provides the estimates of the first time step error. The proof is given in Appendix \ref{sec:Lemma 5}.
\begin{Lem} \label{Lem5}
    Suppose that Hypotheses \ref{Hyp1}, \ref{HypMM}, \ref{Hyp2}, and \ref{Hyp4} hold true. Then, it holds that 
    \begin{equation}
        ||e^1_h||\leq ||e^1_h||+\sqrt{\nu\Delta t}||\nabla e^1_h||\leq C (\Delta t^2 +h^2) ||\phi||_{Z^3\cap H^2(H^2)\cap H^1(H^3)}.\label{Lem5eq}
    \end{equation} 
\end{Lem}
\paragraph{Proof of Theorem \ref{Thm8}}
We substitute $e^n_h$ in the numerical scheme \eqref{LG1} and obtain the following equations for the error:
\begin{align} 
\left(\frac{e^n_h-[e^{n-1}_h\circ X^n_1]\gamma^n}{\Delta t},\psi_h\right)+a_0(e^n_h,\psi_h)&=\langle R_h^n,\psi_h\rangle,  \hspace{0.2cm} \forall\psi_h\in \Psi^n_h, \hspace{0.2cm} n=1, \\
\left( \fz{ 3 e^n_h -4e^{n-1}_h \circ X_1^n \gamma^n + e^{n-2}_h \circ \tilde{X}_1^n \tilde{\gamma}^n}{2\Delta t} ,\psi_h\right)+a_0(e^n_h,\psi_h)&=\langle \tilde{R}^n_h,\psi_h\rangle,  \hspace{0.2cm} \forall\psi_h\in V^n_h, \hspace{0.2cm} n\geq 2, \label{err2}
\end{align}
where the residual $R^n_h,R^n_1,R^n_2,$ and $R^n_3$ are given as in the proof of Theorem \ref{Thm7}, cf., Section \ref{sec:proofslg1}, while the residual on the right hand side of \eqref{err2} is given by:
\begin{align*}
\tilde{R}^n_h &\defeq \tilde{R}^n_1 + \tilde{R}^n_2 + R^n_3, \\
\tilde{R}^n_1 &\defeq \frac{\partial \phi^n}{ \partial t} + \nabla \cdot \ (u^n \phi^n) -  \fz{ 3\phi^n-4\phi^{n-1} \circ X_1^n \gamma^n + \phi^{n-2} \circ \tilde{X}_1^n \tilde{\gamma}^n}{2\Delta t}, \\
\tilde{R}^n_2 &\defeq  \fz{ 3\eta^n-4\eta^{n-1} \circ X_1^n \gamma^n + \eta^{n-2} \circ \tilde{X}_1^n \tilde{\gamma}^n}{2\Delta t}.
\end{align*}
To obtain an estimate for $ \| \tilde{R}_1 \| $, we follow the error estimate framework for the general convection-diffusion problem on uniform mesh (details are given in Appendices~\ref{sec:boundr12} and \ref{sec:boundr3}), which gives us
\begin{align}
 \| \tilde{R}_1 \| _{\ell^2(\Psi'_h)} &\leq  C_1 \Delta t^2  \| \phi \| _{Z^3(0,T)} \label{eq:r1boundlg2}
\end{align}
and as we have shown in Appendix~\ref{sec:boundr3} it holds $R_3^n=0$.
Next, we compute an estimate for $ \| \tilde{R}_2 \| _{(\Psi_h^n)^\prime}$. For $n \geq 2$ it holds
\begin{align*}
     \| \tilde{R}_{2}^n \| _{(\Psi_h^n)^\prime}&=\frac{1}{2\Delta t} \| 3\eta^n - 4\eta^{n-1} \circ X^n_1 \gamma^n + \eta^{n-2} \circ \tilde{X}^n_1 \tilde{\gamma}^n \|_{(\Psi_h^n)^\prime} \notag \\
    &=\left\Vert \frac{3}{2} \bar{D}_{\Delta t } \eta^n - \frac{1}{2}\bar{D}_{\Delta t} \eta^{n-1} + \frac{2}{\Delta t}(\eta^{n-1}-\eta^{n-1}\circ X^n_1\gamma^n) - \frac{1}{2\Delta t}(\eta^{n-2}-\eta^{n-2}\circ \tilde{X}^n_1 \tilde{\gamma}^n) \right\Vert_{(\Psi_h^n)^\prime} \notag \\
    &\leq \frac{3}{2}  \| \bar{D}_{\Delta t}\eta^n \|  + \frac{1}{2}  \| \bar{D}_{\Delta t}\eta^{n-1} \|  +  \frac{2}{\Delta t} \| (\eta^{n-1}-\eta^{n-1}\circ X^n_1\gamma^n) \| _{(\Psi_h^n)^\prime} + \frac{1}{2\Delta t} \| (\eta^{n-2}-\eta^{n-2}\circ \tilde{X}^n_1 \tilde{\gamma}^n) \| _{(\Psi_h^n)^\prime} \notag \\
    &\leq C_3 (  \| \bar{D}_{\Delta t}\eta^n \|  +  \| \bar{D}_{\Delta t}\eta^{n-1} \|  +  \| \eta^{n-1} \|  +  \| \eta^{n-2} \|  ) \quad (\because\text(Lem.\ref{lem:composite})) \notag\\
    &\leq C_4\left(h^2 \Delta t^{-1/2} \| \phi \| _{H^1(t^{n-2},t^n;H^2(\Omega))} +  h^2 \| \phi \| _{H^1(H^3)} \right), 
\end{align*}
where the last inequality follows from Lemma~\ref{lem:interpolation} and Theorem~\ref{Thm2}. Taking the $\ell^2$-norm in the previous estimate we obtain
\begin{align}\label{eq:r2boundlg2}
   \| \tilde{R}_2 \| _{\ell^2(\Psi'_h)} &\leq C_5  h^2 \left( \| \phi \| _{H^1(0,T;H^2(\Omega))} +   \| \phi \| _{H^1(H^3)} \right).
\end{align}
Combining the bounds \eqref{eq:r1boundlg2} and \eqref{eq:r2boundlg2} and taking into account the fact $ R_3 =0$, it follows
\begin{equation}\label{eq:rboundlg2}
 \| \tilde{R}_h \| _{\ell^2(\Psi'_h)}  \leq C ( \Delta t^2  \| \phi \| _{Z^3(0,T)} + h^2  \| \phi \| _{H^1(0,T;H^2(\Omega))} 
 +h^2  \| \phi \| _{H^1(H^3)} ) ,
\end{equation}
where $C>0$ is independent of $h$ and $\Delta t$. 
Finally, from the stability result of Corollary~\ref{Thm4} and using $e^0_h=0$, we get 
\begin{align*}
    ||e_h||_{\ell^\infty(L^2)} + \nu ||\nabla e_h||_{\ell^2(L^2)} &\leq \left(||e^1_h|| + \sqrt{\nu\Delta t} ||\nabla e^1_h||\right) + \left(||e_h||_{\ell_2^\infty(L^2)} + \sqrt{\nu} ||\nabla e^1_h||_{\ell^2_2(L^2)}\right) \\
    & \leq (||e^1_h|| + \sqrt{\nu\Delta t} ||\nabla e^1_h||) +C||\tilde{R}_h||_{\ell^2_2(\Psi'_h)} \\
    & \leq C_1 (\Delta t^2 + h^2) ||\phi||_{Z^3\cap H^2(H^2)\cap H^1(H^3)} 
\end{align*}
which implies the error estimate \eqref{th8}. We employ Lemma \ref{Lem5} for the estimates of the first time step error such that there is no loss of convergence order. \qed

%\\\\\\\\\\\\\\\\\\\\\\\\\\\\\\\\\\\\\\\\\\\\\\\\\\\\\\\\\\\\\\\\\\\
\section{NUMERICAL EXPERIMENTS}
%\\\\\\\\\\\\\\\\\\\\\\\\\\\\\\\\\\\\\\\\\\\\\\\\\\\\\\\\\\\\\\\\\\\
In this section, we present two numerical experiments using the second order LGMM scheme~\eqref{LG2} combined with the moving mesh method \eqref{MM_d} that show the benefits of the new scheme and verify the error estimate from Theorem~\ref{Thm8}. As initial data we take $\phi_h^0=\Pi_h^0 \phi^0$ from the examples below. To compute the integrals that occur in the scheme we employ the quadrature of order nine. Since linear finite element spaces are used in our proposed scheme we do not consider higher order quadrature formulae as proposed e.g., in \cite{bermejo2015modiflagrangaler}. The linear system appearing in \eqref{LG2} and \eqref{MM_d} are iteratively solved using the conjugate gradient (CG) method and successive over-relaxation (SOR) method, respectively. In all experiments we start with an equidistant mesh at the initial time, i.e., for a given $h_0>0$ the points $P^0_1,\dots, P_{N_p}^{0}$ are such that
\begin{equation}\label{eq:uniforminitialmesh}
  P_{j+1}^0 - P_j^0 = h_0, \qquad \forall i \in\{1, \dots, N_p \}.
\end{equation}
The numerical results obtained by the new LGMM scheme are compared to analogous results by the LG scheme with static mesh, which can be interpreted as LGMM scheme with points satisfying $P_j^n=P_j^0$ for all $i\in\{1,\dots, N_p\}$ and $n\in\{1,\dots,N_T\}$ in addition to \eqref{eq:uniforminitialmesh}.
\begin{Ex}\label{ex2}
  We consider the domain $\Omega = (-1,1)$, final time $T=0.5$ and velocity field $u(x,t) = 1+\sin(t-x)$ in problem~\eqref{prob:c-d}. No force field is assumed in this example as we set $f=0$, and for the boundary conditions we set $g=0$. We take the initial value $\phi^0=\phi(\cdot, 0)$ according to the exact solution
\[
  \phi(x,t) = \exp\left(- \frac{1 - \cos (t-x) }{\nu}\right).
  \]
\end{Ex}
We solve Example~\ref{ex2} with diffusion coefficient set to $\nu=0.01$ and $\nu=0.0001$. In the moving mesh method~\eqref{MM_d} we set $\nu_M=\nu$. The integer $N$ determines the discretization of the domain as we choose the initial mesh size $h_0=2/N$. The time increment is linearly coupled to the initial mesh size through the relation $\Delta t = 4 h_0$. 
In this example, since the velocity $u$ does not satisfy Hypothesis \ref{Hyp1}, i.e.,$u_{\vert \Gamma} \neq 0$, the non-overlapping condition (cf. Theorem \ref{Thm5}) might not be met. In this case, we allow the nodal points $\{P_i^n\}_{i=1}^{N_p}$ to extend beyond the domain.

\begin{figure}
  \centering
  \includegraphics{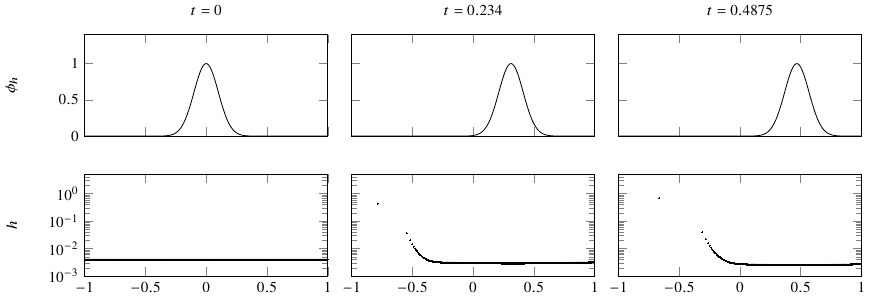}
    \caption{Numerical solution $\phi_h$ and corresponding mesh sizes in Example~\ref{ex2} over the computational domain at time instances $t = 0$~(left), $t=0.2340$~(center) and $t=0.4875$~(right) obtained by the LGMM scheme for $\nu=0.01$ and $N=512$.}\label{fig:numrics2}
  \end{figure}

In Fig.~\ref{fig:numrics2} we show the solution of the LGMM scheme for $N=512$ and $\nu=0.01$ in terms of the functions $\phi^n_h$ together with the corresponding local mesh or partition sizes $h_i^n=P_{i+1}^n-P_i^n$ for $i=1,\dots, N_P$ with respect to their distribution over the computational domain. Clearly, the LGMM scheme maintains a high resolution, i.e., small mesh sizes, in the region, where $\phi_h$ is large, whereas regions with small $\phi_h$ are partly significantly lower resolved. 

Tables~\ref{tab1LG}--\ref{tab2LGMM} show the errors and the corresponding experimental orders of convergence (EOC)\footnote{The EOC is computed by the formula $\text{EOC}=\log_2(E^1/E^2)$ with $E^1$ and $E^2$ denoting the corresponding error in two consecutive lines of the table.} of both the LGMM and the LG scheme of second-order after (initial) grid refinement, i.e., iteratively increasing $N$. In the tables we consider discretization errors with respect to $L^2(\Omega)$, $H^1(\Omega)$ and the loss of total mass, defined as:
\begin{align*}
E_{\ell^\infty(L^2)}\defeq\frac{ \| \phi_h-\Pi_h\phi \| _{\ell^\infty(L^2)}}{ \| \Pi_h\phi \| _{\ell^\infty(L^2)}}, \qquad E_{\ell^2(H^1_0)}&\defeq\frac{ \| \phi_h-\Pi_h\phi \| _{\ell^2(H^1_0)}}{ \| \Pi_h\phi \| _{\ell^2(H^1_0)}}, \quad E_{\rm{mass}}\defeq\frac{\left|\int_\Omega\phi^{N_T}_h dx - \int_\Omega\Pi_h\phi^{N_T}_h dx \right|} {\left|\int_\Omega\Pi_h\phi^{N_T}_h dx\right|},
\end{align*}
where $\| \phi \|_{\ell^2(H_0^1)} \defeq \| \nabla \phi \|_{\ell^2(L^2})$ and $\Pi_h$ denotes the time dependent Lagrange interpolation operator at time instance $t^n$ given as a mapping $\Pi_h(t^n):C(\bar{\Omega})\rarrow\Psi_h^n$. Due to Theorem~\ref{Thm8} and the coupling between $\Delta t$ and $h$ we expected experimental convergence order 2 in both the $\ell^\infty(L^2)$ and the $\ell^2(H_0^1)$ (semi-) norm. While the EOCs in the tables mostly support this expectation a decrease in case of higher mesh resolutions for the LG scheme is visible. In the case $\nu=0.01$ this occurs in $\ell^2(H_0^1)$ and becomes more significant also in $\ell^\infty(L^2)$ in the case $\nu=0.0001$. The LGMM scheme does not suffer from this decrease in EOC and provides in the affected cases more accurate numerical solution in terms of both norms. The tables further exhibit a low relative loss of mass as $E_{\rm{mass}}$ is of low magnitude even for coarse grids and further decreases as the mesh is refined. While the mesh movement of the LGMM scheme leads to slightly larger $E_{\rm{mass}}$ on fine meshes in comparison to the LG scheme if $\nu=0.01$ the loss of mass for the LGMM scheme is significantly smaller than for the LG scheme if $\nu=0.0001$.

\begin{Ex}\label{ex1}
We consider the domain $\Omega = (-1,1)$, final time $T=2$, velocity field $u(x,t) = \sin(2\pi x)$ and diffusion coefficient $\nu = 10^{-5}$ in problem~\eqref{prob:c-d}. Again we take $f=0$ and $g=0$. The initial datum is set to $\phi^0(x) = \exp[- 100(1 - \cos(x))]$.
\end{Ex}

\begin{figure}
  \centering
  \includegraphics{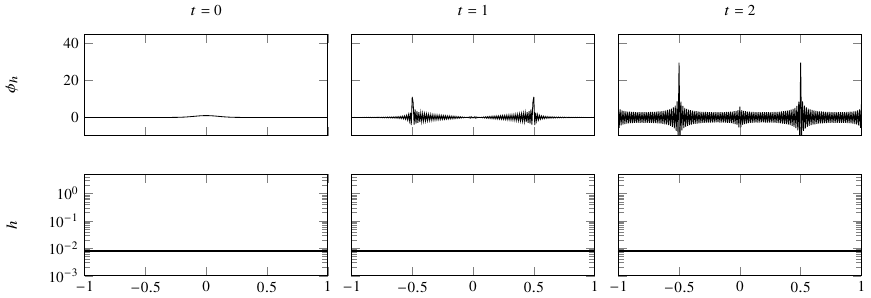}
    \caption{Numerical solution $\phi_h$ and corresponding mesh sizes in Example~\ref{ex1} over the computational domain at time instances $t = 0$~(left), $t=1$~(center) and $t=2$~(right) obtained by the LG scheme with fixed mesh ($N=256$). The numerical solution exhibits oscillations.}\label{fig:numrics1}
  \end{figure}

\begin{figure}
  \centering
  \includegraphics{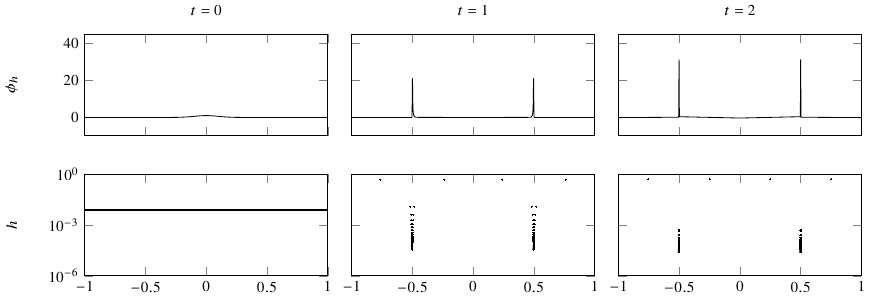}
    \caption{Numerical solution $\phi_h$ and corresponding mesh sizes in Example~\ref{ex1} over the computational domain at time instances $t = 0$~(left), $t=1$~(center) and $t=2$~(right) obtained by the LGMM scheme with $N=256$. The nodal points aggregate along with the solution $\phi_h$.}\label{fig:numrics}
  \end{figure}

We solve Example~\ref{ex1} using scheme~\eqref{LG2} combined with the moving mesh method \eqref{MM_d}, using the parameter $\nu_M=\nu$, an initial uniform mesh satisfying \eqref{eq:uniforminitialmesh} for $h_0=2/1024$ and the fixed time increment~$\Delta t = 10^{-4}$, which satisfies condition \eqref{eq:CFL} during the computation. Again the results by the new LGMM scheme are compared to the LG scheme with static mesh. Comparing Fig.~\ref{fig:numrics1}  and Fig.~\ref{fig:numrics}, we can observe that while the uniform mesh scheme leads to an oscillating solution, the LGMM scheme is capable to capture the high concentration phenomena. This simulation shows the advantage of the proposed LGMM scheme in capturing sharp spike pattern as observed in bio-medical applications. 

\begin{table} 
\centering
\caption{Relative errors and EOCs of LG scheme for $\nu=0.01$}\label{tab1LG}
\begin{tabular}{lllllll}
\hline\noalign{\smallskip}
$N$ & $\Delta t$ & $E_{\ell^\infty(L^2)}$ & EOC & $E_{\ell^2(H^1_0)}$& EOC & $E_{\rm{mass}}$ \\
\noalign{\smallskip}\hline\noalign{\smallskip}
128 & 6.25$\times 10^{-2}$ & 2.795558$\times 10^{-3}$ & - & 4.621785$\times 10^{-3}$ & - & 1.931562$\times 10^{-5}$ \\
256 & 3.12$\times 10^{-2}$ & 8.085728$\times 10^{-4}$ & 1.79 & 1.296162$\times 10^{-3}$ & 1.83 & 1.389046$\times 10^{-6}$\\
512 & 1.56$\times 10^{-2}$ & 2.221100$\times 10^{-4}$ & 1.86 & 3.445636$\times 10^{-4}$ & 1.91 & 1.527363$\times 10^{-6}$\\
1024 & 7.81$\times 10^{-3}$ & 5.927475$\times 10^{-5}$ & 1.91 & 9.098049$\times 10^{-5}$ & 1.92 & 8.199302$\times 10^{-8}$\\
2048 & 3.91$\times 10^{-3}$ & 1.540739$\times 10^{-5}$ & 1.95 & 2.505214$\times 10^{-5}$ & 1.86 & 7.994320$\times 10^{-8}$\\
4096 & 1.95$\times 10^{-3}$ & 3.949271$\times 10^{-6}$ & 1.96 & 7.976085$\times 10^{-6}$ & 1.64 & 8.886055$\times 10^{-8}$\\
\noalign{\smallskip}\hline
\end{tabular}
\end{table}

\begin{table} 
\centering
\caption{Relative errors and EOCs for of LGMM scheme for $\nu=0.01$}\label{tab1LGMM}
\begin{tabular}{lllllll}
\hline\noalign{\smallskip}
$N$ & $\Delta t$ & $E_{\ell^\infty(L^2)}$ & EOC & $E_{\ell^2(H^1_0)}$& EOC & $E_{\rm{mass}}$ \\
\noalign{\smallskip}\hline\noalign{\smallskip}

128 & 6.25$\times 10^{-2}$ & 3.293675$\times 10^{-3}$ & - & 5.441997$\times 10^{-3}$ & - & 1.141478$\times 10^{-6}$ \\
256 & 3.12$\times 10^{-2}$ & 8.756374$\times 10^{-4}$ & 1.91 & 1.467274$\times 10^{-3}$ &  1.87 & 5.715086$\times 10^{-6}$\\
512 & 1.56$\times 10^{-2}$ & 2.265597$\times 10^{-4}$ & 1.95 & 3.853933$\times 10^{-4}$ & 1.93 & 9.131147$\times 10^{-7}$\\
1024 & 7.81$\times 10^{-3}$ & 5.945318$\times 10^{-5}$ & 1.93 & 8.875689$\times 10^{-5}$ & 2.12 & 4.549741$\times 10^{-7}$\\
2048 & 3.91$\times 10^{-3}$ & 1.545287$\times 10^{-5}$ & 1.95 & 2.439410$\times 10^{-5}$ & 1.87 & 4.652903$\times 10^{-8}$\\
4096 & 1.95$\times 10^{-3}$ & 3.948940$\times 10^{-6}$ & 1.96 & 6.051897$\times 10^{-6}$ & 1.98 & 1.034399$\times 10^{-7}$\\
%9.77e-4 & 8.624762$\times 10^{-7}$ & & 4.244640$\times 10^{-6}$ & &5.303902$\times 10^{-8}$
\noalign{\smallskip}\hline
\end{tabular}

\end{table}

\begin{table} 
\centering
\caption{Relative errors and EOCs of LG scheme for $\nu=0.0001$}\label{tab2LG}
\begin{tabular}{lllllll}
\hline\noalign{\smallskip}
$N$ & $\Delta t$ & $E_{\ell^\infty(L^2)}$ & EOC & $E_{\ell^2(H^1_0)}$& EOC & $E_{\rm{mass}}$ \\
\noalign{\smallskip}\hline\noalign{\smallskip}
128 & 6.25$\times 10^{-2}$ & 6.127321$\times 10^{-2}$ & - & 1.255443$\times 10^{-1}$ & - & 1.256237$\times 10^{-2}$ \\
256 & 3.12$\times 10^{-2}$ & 1.369196$\times 10^{-2}$ & 2.16 & 2.916377$\times 10^{-2}$ & 2.10 & 5.927426$\times 10^{-3}$\\
512 & 1.56$\times 10^{-2}$ & 3.286310$\times 10^{-3}$ & 2.06 & 6.026062$\times 10^{-3}$ & 2.27 & 2.317702$\times 10^{-3}$\\
1024 & 7.81$\times 10^{-3}$ & 1.045305$\times 10^{-3}$ & 1.66 & 1.375878$\times 10^{-3}$ & 2.13 & 1.144369$\times 10^{-3}$\\
2048 & 3.91$\times 10^{-3}$ & 5.000259$\times 10^{-4}$ & 1.07 & 5.729551$\times 10^{-4}$ & 1.27 & 5.870192$\times 10^{-4}$\\
4096 & 1.95$\times 10^{-3}$ & 2.650173$\times 10^{-4}$ & 0.91 & 3.289690$\times 10^{-4}$ & 0.78 & 2.986775$\times 10^{-4}$\\
\noalign{\smallskip}\hline
\end{tabular}
\end{table}

\begin{table} 
\centering
\caption{Relative errors and EOCs of LGMM scheme for $\nu=0.0001$}\label{tab2LGMM}
\begin{tabular}{lllllll}
\hline\noalign{\smallskip}
$N$ & $\Delta t$ & $E_{\ell^\infty(L^2)}$ & EOC & $E_{\ell^2(H^1_0)}$& EOC & $E_{\rm{mass}}$ \\
\noalign{\smallskip}\hline\noalign{\smallskip}
128 & 6.25$\times 10^{-2}$ & 1.021001 $\times 10^{-1}$ & - & 2.825924$\times 10^{-1}$ & - & 4.833423 $\times 10^{-4}$ \\
256 & 3.12$\times 10^{-2}$ & 1.898798$\times 10^{-2}$ & 2.42 & 4.633479$\times 10^{-2}$ & 2.60 & 4.112482$\times 10^{-5}$\\
512 & 1.56$\times 10^{-2}$ & 5.634064$\times 10^{-3}$ & 1.75 & 1.141006$\times 10^{-2}$ & 2.02 & 4.714461$\times 10^{-8}$\\
1024 & 7.81$\times 10^{-3}$ & 8.094441$\times 10^{-4}$ & 2.80 & 1.244332$\times 10^{-3}$ & 3.20 & 3.034628$\times 10^{-6}$\\
2048 & 3.91$\times 10^{-3}$ & 2.574393$\times 10^{-4}$ & 1.66 & 6.381969$\times 10^{-4}$ & 0.97 & 5.883933$\times 10^{-7}$\\
4096 & 1.95$\times 10^{-3}$ & 6.442978$\times 10^{-5}$ & 1.99 & 1.598421$\times 10^{-4}$ & 1.99 & 3.556603$\times 10^{-9}$\\
\noalign{\smallskip}\hline
\end{tabular}
\end{table}
%
%
%
%
%
%
%
%
%
%
%\\\\\\\\\\\\\\\\\\\\\\\\\\\\\\\\\\\\\\\\\\\\\\\\\\\\\\\\\\\\\\\\\\\
\section{CONCLUSION}
%\\\\\\\\\\\\\\\\\\\\\\\\\\\\\\\\\\\\\\\\\\\\\\\\\\\\\\\\\\\\\\\\\\\
In this work, we have equipped the mass-preserving Lagrange--Galerkin scheme of second-order in time with a moving mesh method giving rise to the LGMM scheme, which is capable of numerically solving convection-diffusion problems in one space dimension. We also establish the stability and error estimates of the proposed numerical scheme, the latter being with respect to the $\ell^\infty(L^2 )\cap\ell^2(H_0^1)$ -norm, of order $O(\Delta t+h^2)$ if the one-step method is used in time and of order  $O(\Delta t^2+h^2)$ if the two-step scheme is used in time. We show numerical results which support the proved error estimates. The numerical simulations also show that the proposed LGMM scheme is capable to capture high concentration phenomena. In this paper, we have focused on $d=1$, but in principle our scheme is expected to work for the cases~$d=2, 3$. Therefore, in forthcoming research we consider extensions of our scheme to multidimensional problems as well as applications to real-world problems, especially from biology such as immune system dynamics and cancer growth, in which diffusion and aggregation play crucial roles. 

\section*{Acknowledgements}
N.K. thanks the German Science Foundation (DFG) for the financial support through the projects 461365406, 525842915 and 320021702/GRK2326 and further acknowledges the support by the ERS Open Seed Fund of RWTH Aachen University through project OPSF781.
%
%
%
%
%
%
%
%
%
%
%
%

%%%%%%%%%%%%%%%%%%%%%%%%%%%%%%%%%%%%%%%%%%%%%%%%%
%
%
%
%
%
%
%
%
%
%
%
%
% \\\\\\\\\\\\\\\\\\\\\\\\\\\\\\\\\\\\\\\\\\\\\\\\\\\\\\\\\\\\\\\\\\\
\appendix
\section*{APPENDIX}
\section{Proof of Lemma~\ref{Lem3}}\label{sec:proofl1}
%\\\\\\\\\\\\\\\\\\\\\\\\\\
In the following, we assume $t\in[t^{n-1}, t^n]$ and $x\in [P_k(t), P_{k+1}(t)]$ and use the notations $\phi_{k}=\phi(P_{k}(t),t)$ and
$\phi_{k+1}=\phi(P_{k+1}(t),t)$. By our choice of $t$ and $x$ the linear basis functions evaluate as
\begin{equation*}
  \psi_i(x,t)=
  \begin{cases}
    \frac{P_{k+1}(t)-x}{P_{k+1}(t)-P_k(t)} &i=k\\[5pt]
    \frac{x-P_k(t)}{P_{k+1}(t)-P_k(t)} &i=k+1\\[5pt]
    0 & \text{otherwise}
  \end{cases}.
\end{equation*}
Hence, we note the identity
\begin{equation}\label{eq:l1identity}
w^n_k\psi_k +w^n_{k+1}\psi_{k+1}=\frac{x(w^n_{k+1}-w^n_k)+P_{k+1}w^n_k-P_kw^n_{k+1}}{P_{k+1}-P_k}.
\end{equation}
and obtain the time derivatives
\begin{align*}
\frac{\partial}{\partial t}\psi_k
&=\frac{x(P^{'}_{k+1}-P^{'}_k)-P^{'}_{k+1}P_k+P_{k+1}P^{'}_k}{(P_{k+1}-P_k)^2},\\
\frac{\partial}{\partial t}\psi_{k+1}%&=\frac{-P^{'}_k(P_{k+1}-P_k)-(x-P_k)(P^{'}_{k+1}-P^{'}_k)}{(P_{k+1}-P_k)^2},\\
&=-\frac{x(P^{'}_{k+1}-P^{'}_k)+P^{'}_kP_{k+1}-P_kP^{'}_{k+1}}{(P_{k+1}-P_k)^2}.
\end{align*}
Employing the fact that due to \eqref{eq:movingpoints} it holds $P^{\prime}_i(t)=w^n_i$ we compute 
\begin{align*}
  I(x,t)  &=\sum_{i=1}^{N_p}\phi(P_i(t),t) \frac{\partial}{\partial t}\psi_i(x,t) \\
&=\phi_k\frac{\partial}{\partial t}\psi_k +\phi_{k+1}\frac{\partial}{\partial t}\psi_{k+1}\\
          &=\phi_k \frac{x(w_{k+1}^n-w_k^n)-w_{k+1}^n P_k+P_{k+1}w_k^n}{(P_{k+1}-P_k)^2} - \phi_{k+1}\frac{x(w^n_{k+1}-w^n_k)+w_k^n P_{k+1}-P_k w^n_{k+1}}{(P_{k+1}-P_k)^2}\\
          &= - \frac{(\phi_{k+1}-\phi_k)}{(P_{k+1}-P_k)^2}\left[ x (w^n_{k+1}-w^n_k)+(P_{k+1}w^n_k - P_k w^n_{k+1})\right], \\
          &=-\frac{(\phi_{k+1}-\phi_k)}{P_{k+1}-P_k}\left[ w^n_{k+1}\varphi_{k+1}(x,t)+w^n_k\varphi_k(x,t)\right],
\end{align*}
where \eqref{eq:l1identity} has been used in the last step. This concludes the proof of Lemma~\ref{Lem3}. \qed

%\\\\\\\\\\\\\\\\\\\\\\\\\\\\\\\\\\\\\\\\\\\\\\\\\\\\\
\section{Bounds in the $\Psi^\prime$-norm}\label{sec:cor2}
In this appendix we compute estimates for $\|A^n_1\|_{\Psi^\prime}, \|A^n_2\|_{\Psi^\prime}$, and $\|A^n_3\|_{\Psi^\prime}$. To this end we take $v \in H_0^1(\Omega)$ satisfying $\| v \|_{H^1} \leq 1$. We start by showing the estimate \eqref{eq:a1bound} on $\|A^n_1\|_{\Psi^\prime}$ and therefore define for $k\in\{1,\dots,N_p-1\}$ the average $v_k\coloneqq \frac{1}{h_k}\int_{K_k} v(x)\, dx$ and the function
\[
  Q_k(t) \coloneqq \int_0^1  \int_0^{s_0}  (\nabla^2 \phi)(P_k(t)+s_1 h_k(t), t)  ds_1  ds_0
\]
for brevity of notation. Then we estimate
\begin{align*}
     |(A^n_1, v )| 
    &= \biggl| \frac{1}{\Delta t} \int_{t^{n-1}}^{t^n} \sum_{k=1}^{N_p-1} h_k(t) w^n_k  Q_k(t) \,\int_{K_k} \bigl( \psi_k(x,t) - \psi_{k+1}(x,t) \bigr)\, v(x) dx dt \biggr| \\
    &\leq c \frac{h}{\Delta t} \|w\|_{C^0(L^{\infty})}  \int_{t^{n-1}}^{t^n} \| \phi \|_{H^1(H^3)}  \sum_{k=1}^{N_p-1} \int_{K_k} \left| \bigl( \psi_k(x,t) - \psi_{k+1}(x,t) \bigr) \bigl( v(x)-v_k \bigr) \right| dx dt \\
    &\leq  c  \frac{h}{\Delta t} \|w\|_{C^0(L^{\infty})} \| \phi \|_{H^1(H^3)} \int_{t^{n-1}}^{t^n} \sum_{k=1}^{N_p-1} \|\psi_k(\cdot,t) - \psi_{k+1}(\cdot,t)\|_{L^2(K_k)} \|v-v_k\|_{L^2(K_k)} dt\\
    &\leq  c \frac{h}{\Delta t} \|w\|_{C^0(L^{\infty})}  \| \phi \|_{H^1(H^3)}  \int_{t^{n-1}}^{t^n} \left(\sum_{k=1}^{N_p-1} \frac{h_k}{3} \right)^{1/2} \left(\sum_{k=1}^{N_p-1} h_k^2 |v|^2_{H^1(K_k)}\right)^{1/2} dt\\
     &\leq  c' h^2 \|w\|_{C^0(L^{\infty})}  \| \phi \|_{H^1(H^3)} ||v||_{H^1},
\end{align*}
where we have used the Cauchy-Schwartz inequality, the Poincar\'e inequality as well as the following bound and identities
  \begin{align*}
   Q_k(t) \leq c \| \phi \|_{H^1(H^3)}, \quad
    \int_{K_k} (\psi_{k+1}  -\psi_k) \, dx =0, \quad
     \int_{K_k} (\psi_{k+1}  -\psi_k)^2 \, dx = \frac{h_k}{3}. 
  \end{align*} 
Hence, we obtain the estimate
\begin{align*}
  \|A^n_1\|_{\Psi^\prime}= \sup_{\substack{\|v\|_{H^1}\leq 1}} (A^n_1, v ) \leq   c  h^2 
 \|w\|_{C^0(L^{\infty})} \| \phi \|_{H^1(H^3)}. 
\end{align*}

To show the estimate \eqref{eq:a2bound} on $\|A^n_2\|_{\Psi^\prime}$ we note that the bound
\begin{align*}
   \sum_{k=1}^{N_p-1} \int_{K_k} \left| \int_0^1 \int_0^{s_0} \int_{s_1}^{1-s_1}  (\nabla^3 \phi)(P_k(t)+s_2 h_k(t), t)  ds_2 ds_1 ds_0 \, v(x) \right| dx \leq c\|\phi(\cdot,t)\|_{H^3(\Omega)} \|v\|_{L^2(\Omega)}.
\end{align*}
holds and compute
\begin{align*}
  |(A^n_2, v )|  &=  \biggl| \frac{1}{\Delta t}  \int_{t^{n-1}}^{t^n} \sum_{k=1}^{N_p-1} h^2_k(t) \int_{K_k}  \psi_{k+1}(x,t) w^n_k \int_0^1 \int_0^{s_0} \int_{s_1}^{1-s_1}  (\nabla^3 \phi)(P_k(t)+s_2 h_k(t), t)  ds_2 ds_1 ds_0 \, v(x) dx dt \biggr| \\
    &\leq  c \frac{h^2}{\Delta t} \|w\|_{C^0(L^{\infty})}   \int_{t^{n-1}}^{t^n} \|\phi(\cdot,t)\|_{H^3} dt \|v\|_{L^2} \\
    &= c' \frac{h^2}{\sqrt{\Delta t}} \|w\|_{C^0(L^{\infty})} \|\phi\|_{L^2(t^{n-1},t^n; H^3(\Omega))}  \|v\|_{L^2}.
\end{align*}
Embedding $L^2(\Omega)$ in $\Psi^\prime$, we thus obtain the estimate
\begin{align*}
  \|A^n_2\|_{\Psi^\prime} \leq   \frac{c  h^2}{\sqrt{\Delta t}}  \|w\|_{C^0(L^{\infty})} \| \phi \|_{L^2(t^{n-1},t^n; H^3(\Omega))}  \leq  c  h^2  \|w\|_{C^0(L^{\infty})} \| \phi \|_{H^1(H^3)}. 
\end{align*}

Lastly, to show \eqref{eq:a3bound} we define
\[
  R_k(t) \defeq \int_0^1 \int_0^{s_0} (\nabla^2 \phi)(P_k(t) + (1-s_1) h_k(t), t) ds_1 ds_0
\]
and note the bounds
\[
|R_k(t)| \leq c \|\phi\|_{H^1(H^3)}, \qquad
 \left| \int_0^1 (\nabla w^n)(P_k(t) + sh_k(t)) ds \right| \leq c \|w\|_{C^0(W^{1,\infty})},
\]
which allow us to estimate
\begin{align*}
    |(A^n_3, v )| &= \biggl| \frac{1}{\Delta t} \int_{t^{n-1}}^{t^n} \sum_{k=1}^{N_p-1} h^2_k(t) \int_{K_k}   \psi_{k+1}(x,t) \int_0^1 (\nabla w^n)(P_k(t) + sh_k(t)) ds R_k(t) dt \, v(x) \, dx \biggr| \\
    &\leq c  h^2 \|w\|_{C^0(W^{1,\infty})}   \| \phi \|_{H^1(H^3)}  \|v\|_{L^2},
\end{align*}
where we have used the Cauchy--Schwartz inequality. Hence, we obtain
\begin{align*}
  \|A^n_3\|_{\Psi^\prime} \leq   c  h^2  \|w\|_{C^0(W^{1,\infty})} \| \phi \|_{H^1(H^3)},
\end{align*}
which implies \eqref{eq:a3bound}.
\qed

%\\\\\\\\\\\\\\\\\\\\\\\\\\\\\\\\\\\\\\\\\\\\\\\\\\\\
\section{Residual bounds}
In this part of the appendix we derive the bounds for the residual terms that are used in the proofs of Theorems~\ref{Thm7} and \ref{Thm8} in Sections~\ref{sec:proofslg1} and \ref{sec:proofslg2}.
\subsection{Bound for the term $R_1$ in \eqref{err}}\label{sec:boundr11}
To derive an estimate for $R^n_1$ we write $R^n_1=I^n_1+I^n_2$, where
\begin{align*}
    I^n_1 & \defeq \frac{\partial \phi^n}{ \partial t} + u^n \, \nabla\phi^n - \frac{\phi^n - \phi^{n-1}\circ X^n_1}{\Delta t}, \\
    I^n_2 & \defeq  (\nabla \cdot u^n)\phi^n - \frac{\phi^{n-1}\circ X^n_1(1-\gamma^n)}{\Delta t}.
\end{align*}
We first consider the term $I_1^n$. The computation
\begin{align*}
&\phi^n (x) - [\phi^{n-1}\circ X^n_1] (x) \\
&=\phi (x,t^n) - \phi(x-u^n(x) \Delta t, t^{n-1}) \\
&= - \int^1_0 \frac {\partial}{\partial s} [\phi (x-su^n(x) \Delta t, t^n-s\Delta t)] \, ds \\
&=  \int^1_0 u^n(x) \Delta t \, \nabla\phi (x-su^n(x) \Delta t, t^n-s\Delta t) + \Delta t \frac {\partial \phi}{\partial t} (x-su^n(x) \Delta t, t^n-s\Delta t) ds \\
&= \Delta t \int^1_0 \left[\frac {\partial \phi}{\partial t}+u^n(x) \, \nabla\phi \right] (x-su^n(x) \Delta t, t^n-s\Delta t) ds.
\end{align*}
shows that we can write this term in the form
\begin{align*}
I^n_1 (x) &= \int^1_0 \left[\frac{\partial \phi}{\partial t} + u^n (x) \, \nabla\phi\right] (x,t^n) ds - \int^1_0 \left[\frac {\partial \phi}{\partial t}+u^n(x) \, \nabla\phi \right] (x-su^n(x) \Delta t, t^n-s\Delta t) ds \\
&=- \int^1_0 \left.\left[ \left[\frac {\partial \phi}{\partial t}+u^n(x) \, \nabla\phi \right] (x-s_1u^n(x) \Delta t, t^n-s_1\Delta t)\right]\right|^s_{s_1=0} ds \\
&= - \int^1_0 \int^s_0 \frac {\partial}{\partial s_1} \left[\frac {\partial \phi}{\partial t}+u^n(x) \, \nabla\phi \right] (x-s_1u^n(x) \Delta t, t^n-s_1\Delta t) \, ds_1 \, ds \\
&= \Delta t  \int^1_0 \int^s_0 \left(\left[\frac {\partial}{\partial t}+u^n(x) \, \nabla \right]^2\phi\right) (x-s_1u^n(x) \Delta t, t^n-s_1\Delta t) \, ds_1 \, ds.
\end{align*}
Using the Cauchy-Schwartz inequality, we estimate
\begin{align*}
  |I^n_1(x)| &  \leq \Delta t  \int^1_0 \left( \int^s_0 \left(\left[\frac {\partial}{\partial t}+u^n(x) \, \nabla \right]^2\phi\right) (x-s_1u^n(x) \Delta t, t^n-s_1\Delta t)^2 \, ds_1\right)^{1/2}  \left(\int^s_0 ds_1\right)^{1/2} \, ds \\
 &\leq \Delta t  \int^1_0 \left( \int^s_0 \left(\left[\frac {\partial}{\partial t}+u^n(x) \, \nabla \right]^2\phi\right) (x-s_1u^n(x) \Delta t, t^n-s_1\Delta t)^2 \, ds_1\right)^{1/2} \, ds \\
             &\leq \Delta t \left( \int^1_0 \left(\left[\frac {\partial}{\partial t}+u^n(x) \, \nabla \right]^2\phi\right) (x-s_1u^n(x) \Delta t, t^n-s_1\Delta t)^2 \, ds_1\right)^{1/2}.
\end{align*}
Hence, it holds
\begin{align*}
 \| I^n_1(x) \| ^2
&\leq  \Delta t^2 \int_\Omega  \int^1_0 \left(\left[\frac {\partial}{\partial t}+u^n(x) \, \nabla \right]^2\phi\right) (x-s_1u^n(x) \Delta t, t^n-s_1\Delta t)^2 \, ds_1 \, dx \\
&\leq C_1 \Delta t^2 \int_0^1  \int_\Omega \left(\left[\frac {\partial}{\partial t}+\nabla \right]^2\phi\right) (x-s_1u^n(x) \Delta t, t^n-s_1\Delta t)^2 dx \, ds_1.
\end{align*}
Let $y\defeq x-s_1u^n(x)\Delta t$ and $\tau\defeq t^n-s_1\Delta t$. Then, by a change of variable, we have
\begin{align*}
 \| I^n_1(x) \| ^2&\leq C_2 \Delta t^2 \int_0^1  \int_\Omega \left(\left[\frac {\partial}{\partial t}+ \nabla \right]^2\phi\right) (y, t^n-s_1\Delta t)^2 dy \, ds_1 \\
&= C_2 \Delta t^2 \int_{t^{n-1}}^{t^n} \frac{1}{\Delta t}  \int_\Omega \left(\left[\frac {\partial}{\partial t}+ \nabla \right]^2\phi\right) (y, \tau)^2 dy \, d\tau \\
&= C_2\Delta t \left\Vert\left[\frac {\partial}{\partial t}+ \nabla \right]^2\phi \, \right\Vert^2_{L^2(t^{n-1},t^n; L^2(\Omega))} \\
& \leq C_3 \Delta t  \| \phi \| ^2_{Z^2(t^{n-1},t^n)}.
\end{align*}
By taking the square root, we obtain 
\begin{equation*}
 \| I^n_1(x) \| _{L^2(\Omega)} \leq C_4 \sqrt{\Delta t}  \| \phi \| _{Z^2(t^{n-1},t^n)}.
\end{equation*}
On the other hand, we note that $\frac{(1-\gamma^n)}{\Delta t} = \nabla \cdot u^n + \mathcal{O}(\Delta t)$, which using \eqref{eq:compboundL2} leads to
\begin{align*}
     \| I^n_2 \|  &=  \| \nabla \cdot u^n (\phi^n-\phi^{n-1}\circ X^n_1) + [\phi^{n-1}\circ X^n_1] \mathcal{O}(\Delta t ) \| \leq C_5 \sqrt{\Delta t}  \| \phi \| _{Z^1(t^{n-1},t^n)}.
\end{align*}
By combining the estimates of $I^n_1$ and $I^n_2$, we have 
\begin{align*}
 \| R_1 \| _{\ell^2(\Psi'_h)} &\leq \left(\Delta t \sum_{n=1}^{N_T}  \| R^n_1 \| ^2\right)^{1/2}
\leq C_6 \Delta t  \| \phi \| _{Z^2(0,T)},
\end{align*}
where $C_6>0$ is independent of $h$ and $\Delta t$. \qed
%\\\\\\\\\\\\\\\\\\\\\\\\\\\\\\\\\\\\\\\\\\\\\\\\\\\\\\\\\\\\\\\\\\\

\subsection{Bound for the term $\tilde{R}_1$ in \eqref{err2}}\label{sec:boundr12}
First, we note that  $\gamma^n$ and $\tilde{\gamma}^n$ can be written as 
\begin{align*}
    \gamma^n(x)&=1-\Delta t \nabla \, u^n(x), \qquad
    \tilde{\gamma}^n(x)=1-2 \Delta t \nabla \, u^n(x).
\end{align*}
Therefore, the term $\tilde{R}^n_1$ in \eqref{err2} is recasted as
\begin{align*}
    \tilde{R}^n_1&=\frac{1}{2 \Delta t}\left(3\phi^n-4\phi^{n-1}\circ X^n_1\gamma^n + \phi^{n-2} \circ \tilde{X}^n_1 \tilde{\gamma}^n \right) - \left[\frac{\partial \phi}{\partial t} + \nabla \cdot (u^n \phi^n)\right](\cdot,t^n) \\
    &=\left[ \frac{1}{2 \Delta t}\left(3\phi^n-4\phi^{n-1}\circ X^n_1 + \phi^{n-2} \circ \tilde{X}^n_1 \right) - \left(\frac{\partial \phi}{\partial t} + u^n \nabla \cdot \phi^n \right) \right] + \nabla \cdot u^n [2\phi^{n-1}\circ X^n_1 - \phi^{n-2} \circ \tilde{X}^n_1 - \phi^n] \\
    &\eqdef I^n_1 + I^n_2.
\end{align*}
Proceeding in analogy to \ref{sec:boundr12}, let $y(x, s)=y(x,s;n):=x-u^n(x)(1-s)\Delta t$ and $t(s)=t(s;n):=t^{n-1}+s \Delta t$. Then, the terms $I^n_1$ and $I^n_2$ can be expressed in terms of the integrals
\begin{align*}
    I^n_1(x)&=-2\Delta t^2 \int_0^1 s  \int_{2s-1}^s \left[\left(\frac{\partial}{\partial t} + u^n(x) \, \nabla \right)^3 \phi \right] (y(x,s_1),t(s_1))ds_1 \, ds, \\
    I^n_2(x)&=-\Delta t^2 (\nabla \cdot u^n)(x) \int_0^1 \int_{s-1}^s \left[\left(\frac{\partial}{\partial t} + u^n(x) \, \nabla \right)^2 \phi \right] (y(x,s_1),t(s_1))ds_1 ds.
\end{align*}
Now, we can estimate
\begin{align*}
 \| I^n_1 \|  &= 2\Delta t^2 \left\Vert \int_0^1 s  \int_{2s-1}^s \left[\left(\frac{\partial}{\partial t} + u^n(x) \cdot \nabla \right)^3 \phi \right] (y(\cdot,s_1),t(s_1))ds_1 \, ds \right\Vert \\
&\leq C_1\Delta t^2 \int_0^1 s \int_{2s-1}^s \left\Vert\left[\left(\frac{\partial}{\partial t} + \nabla \right)^3 \phi \right] (y(\cdot,s_1),t(s_1)) \right\Vert ds_1 \, ds \\
&\leq C_2\Delta t^2 \int_0^1 s \int_{2s-1}^s \left\Vert\left[\left(\frac{\partial}{\partial t} +  \nabla \right)^3 \phi \right] (\cdot,t(s_1)) \right\Vert ds_1 \, ds \\ 
&\leq C_3 \Delta t \int_{t^{n-2}}^{t^n} \left\Vert\left[\left(\frac{\partial}{\partial t} +  \nabla \right)^3 \phi \right] (\cdot,t) \right\Vert dt \\
&\leq \sqrt{2} C_3 \Delta t^{3/2} \left\Vert\left(\frac{\partial}{\partial t} + \nabla \right)^3 \phi  \right\Vert_{L^2(t^{n-2},t^n;L^2)} \\
&\leq C_4  \Delta t^{3/2}  \| \phi \| _{Z^3(t^{n-2},t^n)}
\end{align*}
and similarly
\begin{align*}
     \| I^n_2 \|  &\leq  C_5\Delta t^2 \int_0^1 s \int_{s-1}^s \left\Vert\left[\left(\frac{\partial}{\partial t} + \nabla \right)^2 \phi \right] (y(\cdot,s_1),t(s_1)) \right\Vert ds_1 \, ds \\
    &\leq C_6 \Delta t \int_{t^{n-2}}^{t^n} \left\Vert\left[\left(\frac{\partial}{\partial t} + 1 \cdot \nabla \right)^2 \phi \right] (\cdot,t) \right\Vert dt \\
    &\leq C_7  \Delta t^{3/2}  \| \phi \| _{Z^2(t^{n-2},t^n)}.
\end{align*}
By combining the bounds of $I^n_1$ and $I^n_2$ and taking the $\ell^2$-norm over all time instances, we obtain
\begin{align*}
     \| \tilde{R}_1 \| _{\ell^2(\Psi'_h)}\leq \left(\Delta t \sum_{n=1}^{N_T} \| R^n_1 \| ^2\right)^{1/2}
    &\leq C_8 \Delta t^{2}  \| \phi \| _{Z^3(0,T)},
\end{align*}
where $C_8>0$ is independent of $h$ and $\Delta t$. \qed
%\\\\\\\\\\\\\\\\\\\\\\\\\\\\\\\\\\\\\\\\\\\\\\\\\\\\\\\\\\\\\\\\\\\
\subsection{Bound for the term $R^n_3$ in \eqref{err} and \eqref{err2}}\label{sec:boundr3}
We compute an estimate for $R^n_3$ in $(\Psi_h^n)^\prime$. To this end let $v_h\in \Psi_h^n$ be such that $\| v_h\|_{H_1} =1$ and further let $K_i=[P_i(t^n), P_{i+1}(t^n)]$ for all $i\in\{1,\dots, N_p-1\}$. Then we have, employing the fact that for $v_h \in \Psi_h^n$ the function $\nabla v_h \rvert_{K_i}$ is constant 
\begin{align*}
 | \langle R^n_3, v_h\rangle | &= | a_0(\eta^n,v_h) |
  \\&\leq |a_0(\phi ,v_h) -  a_0(\Pi_h \phi,v_h)| \\
                         &\leq \nu \sum_{i=1}^{N_p-1} \left| \int_{K_i}  \left(\nabla \phi \nabla v_h - \frac{\phi^n(P_{i+1}) - \phi^n(P_{i}) }{h_i} \nabla v_h \right) \, dx \right|\\
  &= \nu \sum_{i=1}^{N_p-1} |\nabla v_h \rvert_{K_i}| \left| \phi^n(P_{i+1}) - \phi^n(P_{i}) -  \frac{\phi^n(P_{i+1}) - \phi^n(P_{i}) }{h_i} \int_{K_i}\, dx \right| = 0.
\end{align*}
Hence, it follows that $R^n_3=0\in (\Psi_h^n)^\prime$. \qed

\section{Proofs of the mass preserving properties}
To prove Corollaries \ref{Thm1} and \ref{Thm3} we first state a proposition that will be used in the proofs.
\begin{Prop}[ \cite{rui2010massconsercharac,futai2022masspresertwo}]
Suppose that Hypotheses~\ref{Hyp1}, \ref{HypMM} and \ref{Hyp2} hold true. Then it holds that $X^n_1(\Omega)=\tilde{X}^n_1(\Omega)=\Omega$ and $1/2 \leq \gamma^n, \tilde{\gamma}^n\leq 3/2$ for $n=0,\dots,N_T$.
\end{Prop}
%\\\\\\\\\\\\\\\\\\\\\\\\\\\\\\\\\\\\\\\\\\\\\\\\\\\\\\\\\\\\\\\\\\\
\subsection{Proof of Corollary \ref{Thm1}}\label{sec:proofmasslg1}
\paragraph{Proof of Corollary \ref{Thm1}} Suppose that Hypotheses \ref{Hyp1}, \ref{HypMM} and \ref{Hyp2} hold true. Due to Proposition 2, it holds for all $\rho \in \Psi, n=1,\dots,N_T$ that
\begin{equation*}
    \int_\Omega \rho \circ X^n_1(x) \gamma^n(x) dx = \int_\Omega \rho dx.
\end{equation*}
We prove the theorem by induction. Let $m\in\{2,\dots,N_T\}$ and assume that (\ref{th3}) holds true for $n=m-1$. By substituting $1 \in \Psi_h$ into $\psi_h$ in the scheme (\ref{LG}), we obtain
\begin{align*}
    \int_\Omega \phi^m_h(x) dx &= \int_\Omega \phi^0_h \circ X^m_1(x) \gamma^m(x) dx + \Delta t \left(\int_\Omega f^m(x) dx+ \int_\Gamma g^m(x) ds\right) \\
    &= \int_\Omega \phi^0_h(x) dx +  \Delta t \sum_{i=1}^m\left(\int_\Omega f^i(x) dx+ \int_\Gamma g^i(x) ds\right),
\end{align*}
which proves \eqref{th3}. \qed
%\\\\\\\\\\\\\\\\\\\\\\\\\\\\\\\\
\subsection{Proof of Corollary \ref{Thm3}}\label{sec:proofmasslg2}
Suppose that Hypothesis \ref{Hyp1}, \ref{HypMM} and \ref{Hyp2} holds true. By Proposition 2, it holds for all $\rho \in \Psi, n=1,\dots,N_T$ that
\begin{equation*}
    \int_\Omega \rho \circ \tilde{X}^n_1(x) \tilde{\gamma}^n(x) dx = \int_\Omega \rho dx.
\end{equation*}
We prove the corollary by induction. Let $m\in\{2,\dots,N_T\}$ and assume that (\ref{th6}) holds true for $n=m-1$. By substituting $1 \in \Psi_h$ into $\psi_h$ in the scheme (\ref{LG2}) and due to Theorem 3, we obtain
\begin{align*}
    \int_\Omega \left(\frac{3}{2}\phi^m_h - \frac{1}{2}\phi^{m-1}_h\right) dx &=   \int_\Omega \left(\frac{3}{2}\phi^m_h - \frac{1}{2}\phi^{m-1}_h \circ X^m_1 \gamma^m\right) dx \\
    &=   \int_\Omega \left(\frac{3}{2}\phi^{m-1}_h \circ X^m_1 \gamma^m - \frac{1}{2}\phi^{m-2}_h \circ \tilde{X}^m_1 \tilde{\gamma}^m\right) dx +  \Delta t \left(\int_\Omega f^m(x) dx+ \int_\Gamma g^m(x) ds\right) \\
    &=   \int_\Omega \left(\frac{3}{2}\phi^{m-1}_h  - \frac{1}{2}\phi^{m-2}_h \right) dx +  \Delta t \left(\int_\Omega f^m(x) dx+ \int_\Gamma g^m(x) ds\right) \\
    &=   \int_\Omega \phi^0_h(x) dx +  \Delta t \sum_{i=1}^m\left(\int_\Omega f^i(x) dx+ \int_\Gamma g^i(x) ds\right),
\end{align*}
which implies (\ref{th6}). \qed

\section{Proofs of the stability results}
\subsection{Proof of Corollary \ref{Thm2}}\label{sec:proofstabilitylg1}
We substitute $\phi^n_h \in \Psi_h$ into the numerical scheme \eqref{LG} and obtain 
\begin{equation} \label{stb}
\left(\frac{\phi^n_h-(\phi^{n-1}_h\circ X^n_1)\gamma^n}{\Delta t},\phi^n_h \right) + \nu  \| \nabla \phi^n_h \| ^2 = \langle F^n,\phi^n_h \rangle. 
\end{equation}
By Young's inequality, the functional on the right hand side of \eqref{stb} can be estimated as
\begin{equation}\label{eq:stabrhsbound}
\langle F^n,\phi^n_h\rangle \leq \left( \frac{1}{2\nu} \| F^n \| ^2_{(\Psi_h^n)^\prime} + \frac{\nu}{2}  \| \phi^n_h \| ^2_{(\Psi_h^n)^\prime}\right).
\end{equation}
Let the first term on the left hand side of \eqref{stb} be denoted by $I_n$, then the following lower bound holds
\begin{align*}
I_n&\defeq\left(\frac{\phi^n_h-(\phi^{n-1}_h\circ X^n_1)\gamma^n}{\Delta t},\phi^n_h\right) \\
&=\frac{1}{\Delta t} \left[\frac{1}{2} \| \phi^n_h\ \| ^2 - \frac{1}{2}  \| (\phi^{n-1}_h\circ X^n_1)\gamma^n \| ^2 + \frac{1}{2}  \| \phi^n_h-(\phi^{n-1}_h\circ X^n_1)\gamma^n \| ^2  \right] \\
&\geq \frac{1}{\Delta t} \left[\frac{1}{2} \| \phi^n_h\ \| ^2 - \frac{1}{2}  \| (\phi^{n-1}_h\circ X^n_1)\gamma^n \| ^2 \right].
\end{align*}
Since $\gamma^n-1=\mathcal{O}(\Delta t)$ it holds
\begin{align*}
 \| (\phi^{n-1}_h\circ X^n_1)\gamma^n \|  &=  \| (\phi^{n-1}_h\circ X^n_1)\gamma^n -\phi^{n-1}_h\circ X^n_1 + \phi^{n-1}_h\circ X^n_1  \|  \\
&\leq  \| (\phi^{n-1}_h\circ X^n_1)(\gamma^n-1) \|  +  \| (\phi^{n-1}_h\circ X^n_1) \|  \\
&\leq C_1 \Delta t \| (\phi^{n-1}_h\circ X^n_1) \|  +  \| (\phi^{n-1}_h\circ X^n_1) \|  \\
&= (1+C_1 \Delta t)  \| (\phi^{n-1}_h\circ X^n_1) \|  \\
&\leq (1+2 C_2 \Delta t)  \| \phi^{n-1}_h \| .
\end{align*}
Then, $I_n$ can be written as
\begin{equation}\label{eq:stabIlowerbound}
I_n \geq \frac{1}{\Delta t} \left[\frac{1}{2} \| \phi^n_h\ \| ^2 - \frac{1}{2}  \| \phi^{n-1}_h \| ^2 \right] - C_2 \| \phi^{n-1}_h \| ^2 .
\end{equation}
By substituting \eqref{eq:stabIlowerbound} and \eqref{eq:stabrhsbound} into \eqref{stb} we obtain
\begin{equation} \label{stb2}
\frac{1}{\Delta t} \left[\frac{1}{2} \| \phi^n_h\ \| ^2 - \frac{1}{2}  \| \phi^{n-1}_h \| ^2 \right] + \frac{\nu}{2} \| \nabla \phi^n_h \| ^2 \leq  \frac{1}{2\nu} \| F^n \| ^2_{(\Psi_h^n)^\prime} + c" \| \phi^{n-1}_h \| ^2.
\end{equation}
To complete the proof we apply the Gronwall inequality to \eqref{stb2} and get 
\begin{equation*}
 \| \phi_h \| _{\ell^\infty(L^2)} + \sqrt{\nu}  \| \phi_h \| _{\ell^2(H^1_0)} \leq  C \left[ \| \phi^0_h \| + \| F \| _{\ell^2(\Psi'_h)}\right],
\end{equation*}
where the constant $C>0$ is independent of $\Delta t$ and $h$.  \qed
%\\\\\\\\\\\\\\\\\\\\\\\\\\\\\\\\\\\\\\\\\\\\\\\\\\\\\\\\\\\\\\\\\\\
\subsection{Proof of Corollary \ref{Thm4}}\label{sec:proofstabilitylg2}
First we state the following lemma that will be used in the last part of the proof.
\begin{Lem}[Gronwall's Inequality \cite{futai2022masspresertwo}]\label{lem:gronwall}
Let $a_0$, $a_1$ and $a_2$ be non-negative numbers such that $a_1 \geq a_2$, let further $\Delta t \in (0, 3/(4a_0)]$ and $\{x_n\}_{n\geq 0}, \{y_n\}_{n\geq 1}, \{z_n\}_{n\geq 2} $, $\{b_n\}_{n\geq 2}$ be non-negative sequences. Suppose that
\begin{equation} \label{Lem1}
    \frac{1}{\Delta t} \left(\frac{3}{2}x_n - 2x_{n-1} + \frac{1}{2}x_{n-2} + y_n - y_{n-1}\right) + z_n \leq a_ox_n + a_1x_{n-1} + a_2x_{n-2} + b_n, \quad \forall n \leq 2
\end{equation}
is satisfied. Then it holds that
\begin{equation}
    x_n+\frac{2}{3}y_n+\frac{2}{3}\Delta t \sum_{i=2}^{n}z_i \leq (\exp(2(a_0+a_1+a_2)n\Delta t)+1)\left(x_0+\frac{3}{2}x_1+y_1+\Delta t \sum_{i=2}^n b_i \right), \quad \forall n \leq 2.
\end{equation}
\end{Lem} 

\paragraph{Proof of Corollary \ref{Thm4}} For $n \geq 2$, note that the scheme \eqref{eq:lg22} can be written as 
\begin{equation}
\left( \fz{ 3\phi_h^n-4\phi_h^{n-1}  + \phi_h^{n-2}}{2\Delta t} , \psi_h \right) + a_0( \phi_h^n, \psi_h ) = \lA F^n, \psi_h \rA + \lA I^n_h, \psi_h \rA, \quad \forall \psi_h \in \Psi_h^n,
\end{equation}
where $I^n_h \in (\Psi_h^n)^\prime$ is given by
\begin{align*}
    I^n_h &\defeq \frac{1}{2\Delta t} \left[-4(\phi^{n-1}_h - \phi^{n-1}_h \circ X^n_1\gamma^n) + (\phi^{n-2}_h - \phi^{n-2}_h \circ \tilde{X}^n_1 \tilde{\gamma}^n)\right], \\
    &= \frac{1}{2\Delta t} \left[-4(\phi^{n-1}_h - \phi^{n-1}_h \circ X^n_1) + (\phi^{n-2}_h - \phi^{n-2}_h \circ \tilde{X}^n_1 )\right] \\ &+\frac{1}{2\Delta t} \left[-4(\phi^{n-1}_h \circ X^n_1 - \phi^{n-1}_h \circ X^n_1\gamma^n) + (\phi^{n-2}_h \circ \tilde{X}^n_1 - \phi^{n-2}_h \circ \tilde{X}^n_1 \tilde{\gamma}^n)\right] \\
    &\defeq I^n_{h1}+I^n_{h2}
\end{align*}
for $n\in\{2,\dots,N_T\}$.
By substituting $\phi^n_h \in \Psi_h$ as $\psi_h$ into \eqref{LG2} we have 
\begin{equation}  \label{stab}  
\left( \fz{ 3\phi_h^n-4\phi_h^{n-1}  + \phi_h^{n-2} }{2\Delta t} , \phi_h^n \right) + \nu  \| \nabla \phi^n_h \| ^2 = \langle F^n,\phi^n_h \rangle + \langle I^n_h,\phi^n_h \rangle. 
\end{equation}
The first term of the left hand side can thereby be estimated as
\begin{align*}
\left( \fz{ 3\phi_h^n-4\phi_h^{n-1} + \phi_h^{n-2} }{2\Delta t} , \phi_h^n \right) &= \frac{1}{\Delta t}\left[\frac{3}{4} \| \phi^n_h \| ^2 -  \| \phi^{n-1}_h \| ^2 + \frac{1}{4}  \| \phi^{n-2}_h \| ^2 + \frac{1}{4}  \| \phi^{n}_h - 2\phi^{n-1}_h + \phi^{n-2}_h \| ^2 \right] \\
&\quad  + \frac{1}{\Delta t} \left[\frac{1}{2}\left( \| \phi^n_h - \phi^{n-1}_h \| ^2 -  \| \phi^{n-1}_h - \phi^{n-2}_h \| ^2\right)\right] \\
&\geq \frac{1}{\Delta t}\left[\frac{3}{4} \| \phi^n_h \| ^2 -  \| \phi^{n-1}_h \| ^2 + \frac{1}{4}  \| \phi^{n-2}_h \| ^2 \right] \\
&\quad + \frac{1}{\Delta t} \left[ \frac{1}{2}\left( \| \phi^n_h - \phi^{n-1}_h \| ^2 -  \| \phi^{n-1}_h - \phi^{n-2}_h \| ^2\right)\right]. 
\end{align*}
Conversely, the terms on the right hand side can be estimated as
\begin{align*}
\lA F^n,\phi^n_h \rA &\leq  \| F^n \| _{(\Psi_h^n)^\prime} \| \phi^n_h \| _{H^1(\Omega)} \\
&\leq  \| F^n \| _{(\Psi_h^n)^\prime}( \| \phi^n_h) \|  +  \| \nabla \phi^n_h \| ) \\
&\leq \frac{1}{8}  \| \phi^n_h \| ^2 + \frac{\nu}{4}  \| \nabla\phi^n_h \| ^2 + (2+1/\nu) \| F^n \| ^2_{(\Psi_h^n)^\prime}, \\
 \| I^n_{h1} \| _{(\Psi_h^n)^\prime} &\leq C( \| \phi^{n-1}_h +  \| \phi^{n-2}_h \| ) \quad (\because \text{Lemma~\ref{lem:composite} \eqref{eq:compboundH-1}}), \\
 \| I^n_{h2} \|  &\leq \frac{c}{\Delta t} \left( \| \phi^{n-1}_h \circ X^n_1(1-\gamma^n) \|  +  \| \phi^{n-2}_h \circ \tilde{X}^n_1(1-\tilde{\gamma}^n) \| \right) \\
&\leq c_1 ( \| \phi^{n-1}_h \|  +  \| \phi^{n-2}_h \| ) \quad (\because (1-\tilde{\gamma}^n) \leq c_1\Delta t, \text{Lemma~\ref{lem:composite} \eqref{eq:compboundL2nodiff}} ), \\
\lA I^n_h, \phi^n_h \rA &\leq  \| I^n_{h1} \| _{(\Psi_h^n)^\prime} \| \phi^n_h \| _{\Psi_h} +  \| I^n_{h2} \|   \| \phi^n_h \| _{\Psi_h}\\
&\leq  \| I^n_{h1} \| _{(\Psi_h^n)^\prime}( \| \phi^n_h \| + \| \nabla\phi^n_h \| ) +  \| I^n_{h2} \|   \| \phi^n_h \| _{\Psi_h} \\
&\leq (2+\frac{1}{\nu}) \| I^n_{h1} \| _{(\Psi_h^n)^\prime}^2 + 2  \| I^n_{h2} \| ^2 + \frac{1}{4} \| \phi^n_h \| ^2 + \frac{\nu}{4} \| \nabla\phi^n_h \| ^2 \\
&\leq \frac{1}{4} \| \phi^n_h \| ^2 + \frac{\nu}{4} \| \nabla\phi^n_h \| ^2 + c_2 \left(\frac{1}{2} \| \phi^{n-1}_h \| ^2 + \frac{1}{2} \| \phi^{n-2}_h \| ^2 \right).
\end{align*}
Then, by combining the above estimates \eqref{stab} can be rewritten as
\begin{align*}
    &\frac{1}{\Delta t} \left[\frac{3}{4}  \| \phi^n_h \| ^2 -  \| \phi^{n-1}_h \| ^2  + \frac{1}{4}  \| \phi^{n-2}_h \| ^2 + \frac{1}{2} ( \| \phi^n_h - \phi^{n-1}_h \| ^2 -  \| \phi^{n-1}_h - \phi^{n-2}_h \| ^2)\right] + \frac{\nu}{2} \| \nabla \phi^n_h \| ^2 \\
    &\leq \frac {3}{8}  \| \phi^n_h \| ^2 + c_1 \left(\frac{1}{2}  \| \phi^{n-1}_h \| ^2 +  \| \phi^{n-2}_h \| ^2\right) + c_2 \| F^n \| ^2_{(\Psi_h^n)^\prime}.
\end{align*}
To complete the proof of Corollary~\ref{Thm4}, we apply the Gronwall inequality from Lemma \ref{lem:gronwall} and obtain 
\begin{equation*}
      \| \phi_h \| _{\ell^\infty_2(L^2)} + \sqrt{\nu}  \| \nabla\phi_h \| _{\ell^\infty_2(L^2)} \leq C ( \| \phi^0_h \|  +  \| \phi^1_h \|  +  \| F \| _{\ell^2(\Psi'_h)}),
\end{equation*}
where $C>0$ is independent of $\Delta t$ and $h$. By combining this estimate with Corollary~\ref{Thm2}, the proof of Corollary~\ref{Thm4} is completed. \qed

\section{Proof of Lemma \ref{Lem5}} \label{sec:Lemma 5}
Recalling the calculation of $||R^n_1||_{(\Psi^n_h)'}$, cf. \ref{sec:boundr11}, the bound of $||R^n_2||_{(\Psi^n_h)'}$ from (\ref{eq:r2bound_0}), and taking into account the fact $R^n_3=0$, it holds that
    \begin{align}
        ||R^1_h|| &\leq c_1\left(\sqrt{\Delta t}||\phi||_{Z^2(t^0,t^1)}+ \frac{h^2}{\sqrt{\Delta t}}  \| \phi \| _{H^1(t^0,t^1;H^2)} +  h^2 \| \phi \| _{H^1(H^3)}\right) \notag \\
        &\leq c_2 \left(\Delta t ||\phi||_{Z^3}+ h^2 ||\phi||_{H^2(H^2)\cap H^1(H^3)} \right) \leq c_3 (\Delta t + h^2)  ||\phi||_{Z^3\cap H^2(H^2)\cap H^1(H^3)}. \label{lem5_1}
    \end{align}
    By substituting $e^1_h$ into $\psi_h$ in (\ref{err}), dropping the positive term $a_0(e^1_h,e^1_h)$, taking into account $e^0_h=0$, $\langle R^1_h,e^1_h\rangle\leq \|R^1_h\| \|e^1_h\|$, and using (\ref{lem5_1}), we get
     \begin{align*}
        \|e^1_h\| &\leq \Delta t \|R^1_h\| \leq \Delta t c_1 (\Delta t + h^2)  \|\phi\|_{Z^3\cap H^2(H^2)\cap H^1(H^3)} \notag \\
        &\leq c_2(\Delta t^2 + h^2) \|\phi\|_{Z^3\cap H^2(H^2)\cap H^1(H^3)}.
    \end{align*}
    Similarly, substituting $e^1_h$ into $\psi_h$ in (\ref{err}), taking into account $e^0_h=0$, $\langle R^1_h,e^1_h\rangle\leq \|R^1_h\| \|e^1_h\|$, and using (\ref{lem5_1}), we get
    \begin{align*}
        \|e^1_h\|^2 + \nu\Delta t \|\nabla e^1_h\|^2 &\leq \Delta t \|R^1_h\| \|e^1_h\| \\
        &\leq \Delta t c_1 (\Delta t + h^2) (\Delta t^2 + h^2) \|\phi\|^2_{Z^3\cap H^2(H^2)\cap H^1(H^3)} \\
        &\leq c_2(\Delta t^2 + h^2)^2 \|\phi\|^2_{Z^3\cap H^2(H^2)\cap H^1(H^3)},
    \end{align*}
    which implies (\ref{Lem5eq}).\qed

%\\\\\\\\\\\\\\\\\\\\\\\\\\\\\\\\\\\\\\\\\\\\\\\\\\\\\

\begin{thebibliography}{10}
\providecommand{\url}[1]{{#1}}
\providecommand{\urlprefix}{URL }
\expandafter\ifx\csname urlstyle\endcsname\relax
  \providecommand{\doi}[1]{DOI~\discretionary{}{}{}#1}\else
  \providecommand{\doi}{DOI~\discretionary{}{}{}\begingroup
  \urlstyle{rm}\Url}\fi

\bibitem{acharya1990}
Acharya, S., Moukalled, F.H.: An adaptive grid solution procedure for
  convection-diffusion problems.
\newblock J. Comput. Phys. \textbf{91}(1), 32--54 (1990).
\newblock \doi{https://doi.org/10.1016/0021-9991(90)90003-J}

\bibitem{anderson2000mathemmodeltumourinvasmetas}
Anderson, A.R.A., Chaplain, M.A.J., Newman, E.L., Steele, R.J.C., Thompson,
  A.M.: Mathematical {{Modelling}} of {{Tumour Invasion}} and {{Metastasis}}.
\newblock Comput. Math. Method. M. \textbf{2}(2), 129--154 (2000).
\newblock \doi{10.1080/10273660008833042}

\bibitem{baba1981}
Baba, K., Tabata, M.: On a conservative upwind finite element scheme for
  convective diffusion equations.
\newblock RAIRO Anal. Num\'{e}r. \textbf{15}(1), 3--25 (1981)

\bibitem{benitez2012numerlagrangaler}
Ben\'{\i}tez, M., Berm\'{u}dez, A.: Numerical analysis of a second order pure
  {L}agrange-{G}alerkin method for convection-diffusion problems. {P}art {I}:
  {T}ime discretization.
\newblock SIAM J. Numer. Anal. \textbf{50}(2), 858--882 (2012).
\newblock \doi{10.1137/100809982}

\bibitem{benitez2012IInumerlagrangaler}
Ben\'{\i}tez, M., Berm\'{u}dez, A.: Numerical analysis of a second order pure
  {L}agrange-{G}alerkin method for convection-diffusion problems. {P}art {II}:
  {F}ully discretized scheme and numerical results.
\newblock SIAM J. Numer. Anal. \textbf{50}(6), 2824--2844 (2012).
\newblock \doi{10.1137/100809994}

\bibitem{bermejo2012modiflagrangaler}
Bermejo, R., Saavedra, L.: Modified {L}agrange-{G}alerkin methods of first and
  second order in time for convection-diffusion problems.
\newblock Numer. Math. \textbf{120}(4), 601--638 (2012).
\newblock \doi{10.1007/s00211-011-0418-8}

\bibitem{bermejo2015modiflagrangaler}
Bermejo, R., Saavedra, L.: Modified {L}agrange-{G}alerkin methods to integrate
  time dependent incompressible {N}avier-{S}tokes equations.
\newblock SIAM J. Sci. Comput. \textbf{37}(6), B779--B803 (2015).
\newblock \doi{10.1137/140973967}

\bibitem{MR3079104}
Boscheri, W., Dumbser, M.: Arbitrary-{L}agrangian-{E}ulerian one-step {WENO}
  finite volume schemes on unstructured triangular meshes.
\newblock Commun. Comput. Phys. \textbf{14}(5), 1174--1206 (2013)

\bibitem{brenner2008}
Brenner, S.C., Scott, L.R.: The mathematical theory of finite element methods,
  \emph{Texts in Applied Mathematics}, vol.~15, third edn.
\newblock Springer, New York (2008).
\newblock \doi{10.1007/978-0-387-75934-0}

\bibitem{MR2537850}
Carr\'{e}, G., Del~Pino, S., Despr\'{e}s, B., Labourasse, E.: A cell-centered
  {L}agrangian hydrodynamics scheme on general unstructured meshes in arbitrary
  dimension.
\newblock J. Comput. Phys. \textbf{228}(14), 5160--5183 (2009)

\bibitem{carrillo2019kellersegel}
Carrillo, J.A., Kolbe, N., Luk\'{a}\v{c}ov\'{a}-Medvid'ov\'{a}, M.: A hybrid
  mass transport finite element method for {K}eller-{S}egel type systems.
\newblock J. Sci. Comput. \textbf{80}(3), 1777--1804 (2019).
\newblock \doi{10.1007/s10915-019-00997-0}

\bibitem{MR2442406}
Cheng, J., Shu, C.W.: A high order {ENO} conservative {L}agrangian type scheme
  for the compressible {E}uler equations.
\newblock J. Comput. Phys. \textbf{227}(2), 1567--1596 (2007)

\bibitem{chrysafinos2008lagran}
Chrysafinos, K., Walkington, N.J.: Lagrangian and moving mesh methods for the
  convection diffusion equation.
\newblock M2AN Math. Model. Numer. Anal. \textbf{42}(1), 25--55 (2008).
\newblock \doi{10.1051/m2an:2007053}

\bibitem{colera2020lagrangaler}
Colera, M., Carpio, J., Bermejo, R.: A nearly-conservative high-order
  {L}agrange-{G}alerkin method for the resolution of scalar
  convection-dominated equations in non-divergence-free velocity fields.
\newblock Comput. Methods Appl. Mech. Engrg. \textbf{372}, 113366, 24 (2020).
\newblock \doi{10.1016/j.cma.2020.113366}

\bibitem{colera2021lagrangaler}
Colera, M., Carpio, J., Bermejo, R.: A nearly-conservative, high-order, forward
  {L}agrange-{G}alerkin method for the resolution of scalar hyperbolic
  conservation laws.
\newblock Comput. Methods Appl. Mech. Engrg. \textbf{376}, Paper No. 113654, 28
  (2021).
\newblock \doi{10.1016/j.cma.2020.113654}

\bibitem{douglas1982numer}
Douglas Jr., J., Russell, T.F.: Numerical methods for convection-dominated
  diffusion problems based on combining the method of characteristics with
  finite element or finite difference procedures.
\newblock SIAM J. Numer. Anal. \textbf{19}(5), 871--885 (1982).
\newblock \doi{10.1137/0719063}

\bibitem{ewing1981multisgalermethod}
Ewing, R., Russell, T.: Multistep galerkin methods along characteristics for
  convection-diffusion problems.
\newblock In: R.~Vichnevetsky, R.~Stepleman (eds.) Advances in Computer Methods
  for Partial Differential Equations IV, pp. 28--36. IMACS (1981)

\bibitem{futai2022masspresertwo}
Futai, K., Kolbe, N., Notsu, H., Suzuki, T.: A {{Mass-Preserving Two-Step
  Lagrange}}^^e2^^80^^93{{Galerkin Scheme}} for {{Convection-Diffusion
  Problems}}.
\newblock J. Sci. Comput. \textbf{92}(2), 37 (2022).
\newblock \doi{10.1007/s10915-022-01885-w}

\bibitem{gelinas1981}
Gelinas, R.J., Doss, S.K., Miller, K.: The moving finite element method:
  applications to general partial differential equations with multiple large
  gradients.
\newblock J. Comput. Phys. \textbf{40}(1), 202--249 (1981).
\newblock \doi{10.1016/0021-9991(81)90207-2}

\bibitem{huang2011adapt}
Huang, W., Russell, R.D.: Adaptive moving mesh methods, \emph{Applied
  Mathematical Sciences}, vol. 174.
\newblock Springer, New York (2011).
\newblock \doi{10.1007/978-1-4419-7916-2}

\bibitem{hughes1987errat}
Hughes, T.J.R., Franca, L.P., Balestra, M.: Errata: ``{A} new finite element
  formulation for computational fluid dynamics. {V}. {C}ircumventing the
  {B}abu\v{s}ka-{B}rezzi condition: a stable {P}etrov-{G}alerkin formulation of
  the {S}tokes problem accommodating equal-order interpolations''.
\newblock Comput. Methods Appl. Mech. Engrg. \textbf{62}(1), 111 (1987).
\newblock \doi{10.1016/0045-7825(87)90092-2}

\bibitem{jimack1991tempor}
Jimack, P.K., Wathen, A.J.: Temporal derivatives in the finite-element method
  on continuously deforming grids.
\newblock SIAM J. Numer. Anal. \textbf{28}(4), 990--1003 (1991).
\newblock \doi{10.1137/0728052}

\bibitem{kolbe2022adaptrectanmesh}
Kolbe, N., Sfakianakis, N.: An adaptive rectangular mesh administration and
  refinement technique with application in cancer invasion models.
\newblock J. Comput. Appl. Math. \textbf{416}, Paper No. 114442, 18 (2022).
\newblock \doi{10.1016/j.cam.2022.114442}

\bibitem{lukacova-medvidova2017numeroseenpeter}
Luk\'{a}\v{c}ov\'{a}-Medvi\v{d}ov\'{a}, M., Mizerov\'{a}, H., Notsu, H.,
  Tabata, M.: Numerical analysis of the {O}seen-type {P}eterlin viscoelastic
  model by the stabilized {L}agrange-{G}alerkin method. {P}art {I}: {A}
  nonlinear scheme.
\newblock ESAIM Math. Model. Numer. Anal. \textbf{51}(5), 1637--1661 (2017).
\newblock \doi{10.1051/m2an/2016078}

\bibitem{lukacova-medvidova2017numeroseenpeterII}
Luk\'{a}\v{c}ov\'{a}-Medvi\v{d}ov\'{a}, M., Mizerov\'{a}, H., Notsu, H.,
  Tabata, M.: Numerical analysis of the {O}seen-type {P}eterlin viscoelastic
  model by the stabilized {L}agrange-{G}alerkin method. {P}art {II}: {A} linear
  scheme.
\newblock ESAIM Math. Model. Numer. Anal. \textbf{51}(5), 1663--1689 (2017).
\newblock \doi{10.1051/m2an/2017032}

\bibitem{mcrae1982numer}
McRae, G.J., Goodin, W.R., Seinfeld, J.H.: Numerical solution of the
  atmospheric diffusion equation for chemically reacting flows.
\newblock J. Comput. Phys. \textbf{45}(1), 1--42 (1982).
\newblock \doi{10.1016/0021-9991(82)90101-2}

\bibitem{notsu2016errorestimstabil}
Notsu, H., Tabata, M.: Error estimates of a stabilized {Lagrange}--{Galerkin}
  scheme for the {Navier}--{Stokes} equations.
\newblock ESAIM Math. Model. Numer. Anal. \textbf{50}(2), 361--380 (2016).
\newblock \doi{10.1051/m2an/2015047}

\bibitem{pironneau1989finitelemenmethodfluid}
Pironneau, O.: Finite Element Methods for Fluids.
\newblock John Wiley and Sons, Chichester (1989)

\bibitem{rui2002}
Rui, H., Tabata, M.: A second order characteristic finite element scheme for
  convection-diffusion problems.
\newblock Numer. Math. \textbf{92}(1), 161--177 (2002).
\newblock \doi{10.1007/s002110100364}

\bibitem{rui2010massconsercharac}
Rui, H., Tabata, M.: A {{Mass-Conservative Characteristic Finite Element
  Scheme}} for {{Convection-Diffusion Problems}}.
\newblock J. Sci. Comput. \textbf{43}(3), 416--432 (2010).
\newblock \doi{10.1007/s10915-009-9283-3}

\bibitem{sulman2011optim}
Sulman, M., Williams, J.F., Russell, R.D.: Optimal mass transport for higher
  dimensional adaptive grid generation.
\newblock J. Comput. Phys. \textbf{230}(9), 3302--3330 (2011).
\newblock \doi{10.1016/j.jcp.2011.01.025}

\bibitem{tabata1977}
Tabata, M.: A finite element approximation corresponding to the upwind finite
  differencing.
\newblock Mem. Numer. Math. \textbf{4}, 47--63 (1977)

\bibitem{tabata2016lagrangaler}
Tabata, M., Uchiumi, S.: A genuinely stable {L}agrange-{G}alerkin scheme for
  convection-diffusion problems.
\newblock Jpn. J. Ind. Appl. Math. \textbf{33}(1), 121--143 (2016).
\newblock \doi{10.1007/s13160-015-0196-2}

\bibitem{TAKASHI1992115}
Takashi, N., Hughes, T.J.: An arbitrary lagrangian-eulerian finite element
  method for interaction of fluid and a rigid body.
\newblock Comput. Methods Appl. Mech. Eng. \textbf{95}(1), 115--138 (1992)

\bibitem{MR3030556}
\v{C}esenek, J., Feistauer, M., Hor\'{a}\v{c}ek, J., Ku\v{c}era, V.,
  Prokopov\'{a}, J.: Simulation of compressible viscous flow in time-dependent
  domains.
\newblock Appl. Math. Comput. \textbf{219}(13), 7139--7150 (2013)

\end{thebibliography}
\end{document}